\documentclass[a4paper,12pt]{amsart}
\usepackage{amssymb}
\usepackage{xy}
\xyoption{all}
\xyoption{poly}

\newtheorem{theo}{Theorem}[subsection]
\newtheorem{lemma}[theo]{Lemma}
\newtheorem{coro}[theo]{Corollary}
\newtheorem{prop}[theo]{Proposition}
\newtheorem{fact}[theo]{Fact}
\newtheorem{defi}[theo]{Definition}

\theoremstyle{definition}

\newtheorem{question}[theo]{Question}

\def\dim{\operatorname{dim}\nolimits}

\def\parag{{\S}}

\def\CC{{\mathcal{C}}}
\def\CS{{\mathcal{S}}}

\def\CL{{\mathcal{L}}}

\def\CR{{\mathcal{R}}}

\def\CA{{\mathcal{A}}}

\def\CR{{\mathcal{R}}}
\def\CL{{\mathcal{L}}}

\def\CS{{\mathcal{S}}}

\def\BB{{\mathbf{B}}}
\def\BC{{\mathbf{C}}}

\def\BS{{\mathbf{S}}}

\def\Bs{{\mathbf{s}}}

\def\BS{{\mathbf{S}}}

\def\BC{{\mathbb{C}}}
\def\BR{{\mathbb{R}}}
\def\BZ{{\mathbb{Z}}}

\def\BT{{\mathbf{T}}}
\def\Bs{{\mathbf{s}}}
\def\Bt{{\mathbf{t}}}
\def\bt{{\mathbf{t}}}

\def\BM{{\mathbf{M}}}
\def\BG{{\mathbf{G}}}

\def\bs{{\mathbf{s}}}
\def\bu{{\mathbf{u}}}
\def\bv{{\mathbf{v}}}
\def\bw{{\mathbf{w}}}

\def\bw{{\mathbf{w}}}
\def\bc{{\mathbf{c}}}

\def\op{\operatorname{op}\nolimits}
\def\Id{\operatorname{Id}\nolimits}

\def\GL{\operatorname{GL}\nolimits}

\def\Red{\operatorname{Red}\nolimits}

\def\Vis{{\mathcal{V}}}

\def\reg{\operatorname{reg}\nolimits}

\def\im{\operatorname{im}\nolimits}
\def\re{\operatorname{re}\nolimits}

\def\codim{\operatorname{codim}\nolimits}
\def\ie{{\em i.e.}}

\title{The dual braid monoid}
\author{David Bessis}
\address{David Bessis, Institut Girard Desargues,
CNRS UMR 5028 / Universit\'e Lyon 1, 21 avenue Claude Bernard,
69622 Villeurbanne Cedex, France}
\email{bessis@desargues.univ-lyon1.fr}
\begin{document}
\begin{abstract}
We construct a new monoid structure for 
Artin groups associated with finite Coxeter systems.
This monoid shares with the classical positive braid monoid
a crucial algebraic property: it is a Garside monoid.
The analogy with the classical construction indicates there
is a ``dual'' way of studying Coxeter systems, where the pair
$(W,S)$ is replaced by $(W,T)$, with $T$ the set of all reflections.
In the type $A$ case, we recover the monoid constructed by Birman-Ko-Lee.
\end{abstract}
\maketitle    

{\flushleft \bf Introduction}

Combinatorics of Coxeter systems provide
very powerful tools to understand finite real reflection groups,
their geometry and their braid groups. The goal of this article is to describe
an alternate approach to finite real reflection groups. This approach
can be seen as a natural ``twin'' or ``dual'' of the classical theory
of Coxeter groups and Artin groups. The starting point
is to naively replace a finite Coxeter system $(W,S)$ by the pair $(W,T)$,
where $T$ is the set of {\em all} reflections in $W$. Surprisingly, many 
crucial properties of $(W,S)$ can be replaced by analog properties for
$(W,T)$: there is a corresponding natural presentation for $W$ and for
its braid group, there is a good notion of positive braid monoid, which
injects in the braid group and which enjoys a nice normal form; this
gives a new automatic structure for the braid group, a new solution to
the word problem, a coherence rule for actions on categories, finite 
simplicial $K(\pi,1)$'s... The interest is
that although all these objects are algebraic analogs of the classical
ones, they are not isomorphic to them. The geometry of real reflection
groups has two (or more?) nice ways to be studied. The classical approach,
with walls and chambers, put the emphasis on the real structure of hyperplane
arrangement, while the dual combinatorics encode the symmetries of the
complexified arrangement, when looked at from an eigenvector for a Coxeter
element.

Let $(W,S)$ be a finite Coxeter system. The group $W$ is given
by the group presentation
$$< S | \forall s\in S, s^2 =1 \; ; \; \forall s,t \in S,
\underbrace{sts\dots}_{m_{s,t}} = \underbrace{tst\dots}_{m_{s,t}}
>_{\text{group}}.$$
Let $\BB(W,S)$ be the corresponding Artin group. To have simple yet precise
notations, it is convenient to introduce a formal copy $\BS\simeq S$. For
each $\bs\in \BS$, we write $s$ the corresponding element of $S$.
With this convention, $\BB(W,S)$ is defined as the abstract group
$$\BB(W,S) := 
< \BS | \forall \Bs,\Bt \in \BS, \underbrace{\Bs\Bt\Bs\dots}_{m_{s,t}}
= \underbrace{\Bt\Bs\Bt\dots}_{m_{s,t}} >_{\text{group}}.$$
The map $\bs\mapsto s$ extends to a surjective morphism $p:\BB(W,S)
\rightarrow W$.

Since the defining relations are between positive words, the presentation 
of $\BB(W,S)$ can
also be seen as a monoid presentation. We set
$$\BB_+(W,S):=
< \BS | \forall \Bs,\Bt \in \BS, \underbrace{\Bs\Bt\Bs\dots}_{m_{s,t}}
= \underbrace{\Bt\Bs\Bt\dots}_{m_{s,t}} >_{\text{monoid}}.$$
This monoid is often called the {\em positive braid monoid}. We prefer here
the term of {\em classical braid monoid} (short for
Artin-Brieskorn-Deligne-Garside-Saito-Tits monoid).

The structure of $\BB(W,S)$ and $\BB_+(W,S)$ 
has been studied in great detail in 
\cite{deligne} and \cite{brieskornsaito}. 
One of the main
results is that $\BB_+(W,S)$ satisfies the {\em embedding property},
\ie, the morphism $\BB_+(W,S) \rightarrow \BB(W,S)$
is injective. 
In other words, $\BB_+(W,S)$ is isomorphic to the
submonoid of $\BB(W,S)$ generated by $\BS$. This explains both the
notation $\BB_+(W,S)$ and why we did not bother to introduce
 another formal copy of $\BS$ when defining the classical braid monoid.
Another important result is the existence of a nice normal form in
$\BB(W,S)$, which, for example, gives practical solutions to the word problem.

We introduce a new positive presentation for $\BB(W,S)$.
The corresponding monoid, the {\em dual braid monoid},
also satisfies the embedding property and, more generally,
shares with the classical monoid its crucial structural properties.
When $(W,S)$ is of type $A$, we recover the monoid discovered
by Birman-Ko-Lee. 

The correct algebraic setting to handle both the classical monoid and the
dual monoid is the notion of Garside monoid introduced by Dehornoy-Paris.
Both the classical monoid and the dual monoid are Garside monoids, from which
their other properties follow.

The structure of this article is as follows:

In the first section, we study the combinatorics 
of the pair $(W,T)$. We obtain an analog of the Matsumoto property.

In the second section, we start by defining {\em dual braid relations}.
The dual braid monoid is the monoid defined by these (positive)
relations, or equivalently by a {\em pre-monoid} constructed from $(W,T)$
(Theorem \ref{ruth}). We show that the group defined
by the dual braid relations is again $\BB(W,S)$ (Theorem \ref{mainA}).
The dual braid monoid is a Garside monoid (Theorem \ref{mainB}) and 
the whole theory is an algebraic analog of the classical one.
Unfortunately, some of our proofs are still case-by-case (the
exceptional types are done by computer).

The third section is a geometric interpretation. To each regular vector $v$
in a complexified 
realization of $W$ as a reflection group, we associate a certain
{\em local monoid} $M_v$ encoding the affine geometry of the
hyperplane arrangement, seen from $v$. The structure of $M_v$ depends
on the position of $v$ with respect to a certain
stratification. 
When $v$ is chosen in a real chamber, $M_v$ is isomorphic
to $\BB_+(W,S)$. Another particular stratum yields the dual braid monoid.

The geometry and combinatorics arising in types $A$, $B$ and $D$ are
made explicit in the fourth section. 

The last two sections contain complements and applications.

In an appendix, we give a survey of some techniques and results
from the theory of Garside monoids.

\medskip

{\bf \flushleft Note.}
After the first version of the present paper was circulated,
the author was informed by
T. Brady and C. Watt that they were working on the same problem.
They have independently obtained some of our results,
namely the lattice structure of $P_c$ when $W$ is of type
$B$ or $D$ (see Theorem \ref{mainB} below) as well as explicit embeddings
of the monoids in the corresponding braid groups (\cite{bradywatt});
The types $I_2$ and $H_3$
had also been studied independently by Brady. 
\medskip

{\flushleft \bf Acknowledgements.}
The present work builds on a previous collaboration with
Jean Michel and Fran\c cois Digne. Theorem \ref{ruth} was a 
conjecture suggested by Ruth Corran. The author thanks the
three of them for many discussions, suggestions and comments
which led to improvements of this text.

\section{Reduced $T$-decompositions}

This section contains the first steps of what could be a ``dual Coxeter
theory''.

\subsection{Reflection groups}
We call {\sf abstract (finite real) reflection group} a pair $(W,T)$ where
$W$ is a finite group, $T$ a generating subset of $W$ and
there exists a faithful representation $\rho:W \hookrightarrow \GL(V_{\BR})$, 
with $V_{\BR}$ a finite dimensional $\BR$-vector space, satisfying
$$\forall w\in W, \codim( \ker(\rho(w)-\Id)) =1 \Leftrightarrow w\in T.$$
The group $\rho(W)$ is a {\sf geometric (finite real)
reflection group}, with set of reflections $\rho(T)$. We say that
$\rho$ is a {\sf realization} of $W$.

Unless otherwise specified, all reflection groups considered in this
paper are finite and real, and all Coxeter systems are finite 
(``spherical type'').

Since geometric reflection groups are classified by (finite)
Coxeter systems,
all abstract reflection groups can be obtained as follows:
let $(W,S)$ be a (finite) Coxeter system; let $T$ be the closure of $S$
under conjugation; then $(W,T)$ is an abstract reflection group.
Conversely, if $(W,T)$ is an abstract reflection group, one may
always choose $S\subset T$ such that $(W,S)$ is a Coxeter system.
The type of $(W,S)$ does not depend on the choice of $S\subset T$.
The {\sf rank} of $(W,T)$ is the rank $|S|$ of $(W,S)$.

\begin{question}
Is there a nice combinatorial description of abstract reflection groups,
similar to Coxeter systems, allowing for example a direct classification
(not using the classification of Coxeter systems)?
\end{question}

We do not have an answer to this question, but we do obtain here some strong
combinatorial properties of $(W,T)$.

\subsection{The reflection length $l_T$}
An (abstract) reflection group $(W,T)$ is a particular example of
{\em generated group},
as this notion is defined in paragraph \ref{ggroup} of the Appendix.
We have a notion of reduced $T$-decomposition,
a length function $l_T$ and two partial orders $\prec_T$ and $\succ_T$ on
$W$ (see \ref{ggroup}). The function $l_T$ is called
{\sf reflection length}.
Since $T$ is invariant by conjugation, it is
clear that $\prec_T$ and $\succ_T$ coincide.

Carter gave a geometric interpretation of the function $l_T$:
\begin{lemma}
\label{carter}
Let $\rho$ be a realization of a reflection group $(W,T)$.
\begin{itemize}
\item[(i)] Let $w\in W$ and $t\in T$. We have
$$t \prec_T w \Leftrightarrow \ker(\rho(t) -\Id) \supset \ker(\rho(w)-\Id)$$
\item[(ii)] For all $w\in W$, $l(w)=\codim(\ker(\rho(w)-\Id))$.
\end{itemize}
\end{lemma}

\begin{proof}
See \cite{carter}, Lemma 2.8. (Carter actually works
with Weyl groups, but his argument can be used with an arbitrary finite
geometric reflection group).
\end{proof}

\subsection{Chromatic pairs and Coxeter elements}
\begin{defi} A {\sf chromatic pair} for an (abstract) reflection group
$(W,T)$
is an ordered pair $(L,R)$
of subsets of $T$, such that:
\begin{itemize}
\item the intersection $L\cap R$ is empty;
\item the subgroups $<L>$ and $<R>$ are abelian;
\item the pair $(W,L\cup R)$ is a Coxeter system.
\end{itemize}
\end{defi}

When unambiguous, we will sometimes write the pair $L\cup R$
instead of $(L,R)$.
The term ``chromatic'' comes from the fact that
the Coxeter graph of $(W,L\cup R)$ comes equipped with a $2$-colouring:
elements of $L$ are said to be ``left'' (let us pretend this is a colour), 
elements of $R$ are ``right''.
If $(W,S)$ is an irreducible Coxeter system, there are exactly 
two $2$-colourings of the Coxeter graph of $(W,S)$.

If $L\cup R$ is a chromatic pair, we set
$$s_L:= \prod_{s\in L} s, \quad 
s_R:= \prod_{s\in R} s, \quad c_{L,R} := s_Ls_R.$$

\begin{defi}
The {\sf Coxeter elements} of $(W,T)$ are the elements of the form
$c_{L,R}$, where $(L,R)$ is a chromatic pair.
A {\sf dual Coxeter system} is a triple $(W,T,c)$ where $(W,T)$ is a
reflection group, and $c$ is a Coxeter
element in $(W,T)$.
\end{defi}

Coxeter elements form a single conjugacy class. The order of Coxeter 
elements is the {\sf Coxeter number}, denoted by $h$.
Contrary to \cite{bourbaki}, our definition is not specific to a choice
of $S$ -- what we call a Coxeter element is, in the terminology of
\cite{bourbaki}, an ``element of the Coxeter class''.

In the ``dual'' approach, choosing a Coxeter element $c$ plays a similar
role as choosing a Coxeter generating set $S$ (or, in geometric terms, 
a chamber) in the classical approach.

\begin{lemma}
Let $(W,T,c)$ be a dual Coxeter system of rank $n$.
We have $l_T(c)=n$ and $\forall t\in T, t\prec_T c$.
\end{lemma}

\begin{proof}
Let $\rho$ be an essential realization of $W$.
We have $\ker(c-\Id) = \{1\}$ (this is a consequence of
\cite{bourbaki}, Ch. V, \S 6, Th. 1, p. 119).
The result then follows from  Lemma \ref{carter}.
\end{proof}

The last statement of the above lemma will be refined in
\ref{debut}.

\begin{lemma}
\label{orbit}
Let $(W,T)$ be an irreducible reflection group, with Coxeter number $h$.
Let $(L,R)$ be a chromatic pair, let $S:=L\cup R$. 
Then $T$ is the closure of $S$ under the conjugacy action of $c_{L,R}$.
Moreover, if $\Omega \subset T$ is an orbit for the conjugacy action of
$c_{L,R}$, then:
\begin{itemize}
\item[(i)] either $\Omega$ has cardinal $h$ and $\Omega \cap S$ has
cardinal $2$;
\item[(ii)] or $\Omega$ has cardinal $h/2$ and $\Omega \cap S$ has
cardinal $1$.
\end{itemize}
\end{lemma}

\begin{proof}
Write $L=\{s_1,\dots,s_k\}$, $R=\{s_{k+1},\dots,s_n\}$, and
$c:=c_{L,R}=s_1\dots s_n$.

Let $s_i,s_j\in S$.
Assume $s_ic^m=c^ms_j$, for some integer $m>0$. Then
we have $m\geq \lfloor h/2 \rfloor$.

Indeed, assume that $m< \lfloor h/2 \rfloor$; we will find
a contradiction. According to
\cite{bourbaki}, Ch. V, \S 6, Ex. 2 (p 140), 
$$(s_1\dots s_n)^{m}$$ is a reduced $S$-decomposition of $c^{m}$,
and $$s_{k+1}\dots s_n (s_1\dots s_n)^m s_1\dots s_k$$ is 
a reduced $S$-decomposition of $c_{R,L}^{m+1}$.
Assume that $s_i\in L$ ($i\leq k$). Then
$$s_1\dots \hat{s}_i \dots s_n (s_1\dots s_n)^{m-1}$$ and  
$$(s_1\dots s_n)^{m-1} s_1\dots \hat{s}_j \dots s_n$$ are reduced
$S$-decompositions of $s_ic^m=c^ms_j$. But $s_i$ is a left descent of 
only one of these two words, which gives a contradiction. 
Now assume that $s_i$ is right.
Then $s_i(s_1\dots s_n)^m$ is $S$-reduced. 
Since $s_i(s_1\dots s_n)^m = (s_1\dots s_n)^m s_j$, the word
$(s_1\dots s_n)^m s_j$ is also $S$-reduced, and $s_j\in L$.
The word $s_i(s_1\dots s_n)^m s_j$ must also be $S$-reduced 
(view it as a subword of the reduced $S$-decomposition
of $c_{R,L}^{m+1}$ given above).
Since $s_i(s_1\dots s_n)^m s_j=(s_1\dots s_n)^m$, we have a contradiction.

This implies that for each orbit $\Omega$, we have
$|\Omega| \geq h/2 |\Omega \cap S|$. Since $c$ has order $h$, we
also have $|\Omega| \leq h$. Using the well-known relation
$|T| = h/2 |S|$, we obtain the claimed results.
\end{proof}

\subsection{Parabolic Coxeter elements}

\begin{defi}
Let $(W,T)$ be a reflection group. Let $S\subset T$ be such
that $(W,S)$ is a Coxeter system. Let $I\subset S$. Let $W_I:=<I>$
and $T_I:=T\cap W_I$. The reflection group $(W_I,T_I)$ is
a {\sf parabolic subgroup}
of $(W,T)$.

An element $w\in W$ is
a {\sf parabolic Coxeter element} if it is a Coxeter element
in some parabolic subgroup of $(W,T)$.
\end{defi}

\begin{lemma}
\label{debut}
Let $(W,T,c)$ be a dual Coxeter system of rank $n$.
Let $t\in T$. 
There exists a chromatic pair $(L,R)$ such that $t\in L$ and
$c=c_{L,R}$. In particular,
there exists $(t_1,\dots,t_n)\in \Red_T(c)$ such that
$t_1=t$ and $(W,\{ t_1,\dots,t_n\})$ is a Coxeter system.
\end{lemma}

\begin{proof}
Let $(L,R)$ be a chromatic pair such that $c=c_{L,R}$.
Let $t\in T$.
By Lemma \ref{orbit}, $t$ is of the form $c^k s c^{-k}$, with $s\in L\cup R$.
\begin{itemize}
\item Assume $s\in L$. Then $(L',R'):=(c^kLc^{-k},c^kRc^{-k})$ is as
required.
\item If $s\in R$, we note that $c= c_{L,R} = c_{R, s_R^{-1}Ls_R}$,
so, by modifying the chromatic pair, we are back to the case
already discussed.
\end{itemize}
\end{proof}

This allows the following characterization of parabolic Coxeter elements.

\begin{lemma}
\label{paracox}
Let $(W,T)$ be a dual Coxeter system. Let $w\in W$.
The following assertions are equivalent:
\begin{itemize}
\item[(i)] There exists a Coxeter element $c\in W$, such that
$w \prec_T c$.
\item[(ii)] The element $w$ is a parabolic Coxeter element.
\end{itemize}
\end{lemma}

\begin{proof}
Let $c$ be a Coxeter element in $W$. Let $t\in T$.
By the previous lemma, we can find $(t_1,\dots,t_n)\in \Red_T(c)$ such
that $t_1=t$ and $(W,\{t_1,\dots,t_n\})$ is a Coxeter system.
Thus $tc=t_2\dots t_n$ is a Coxeter element in the parabolic subgroup
generated by $\{ t_2,\dots,t_n\}$. By induction, this proves
$(i)\Rightarrow (ii)$. The converse is easy.
\end{proof}

\subsection{Generating sets closed under conjugation}

Let $(G,A)$ be a generated group and assume that
$A$ is invariant by conjugation.
Let $n$ be a positive integer. Consider the Artin group $B_n$ 
of type $A_{n-1}$:
$$B_n \; : \quad \xy  (10,0) ;
 (20,0)  **@{-}; (30,0) *++={\dots} **@{-} ;
 (40,0) **@{-};
    (10,-5) *++={\sigma_1} ; (20,-5) *++={\sigma_2};
    (40,-5) *++={\sigma_{n-1}},
     (10,0) *++={\bullet} ; (20,0) *++={\bullet} ;
      (40,0) *++={\bullet} 
    \endxy .$$
It is clear that the assignement
$$\mbox{\small $\sigma_i (a_1,\dots,a_{i-1},a_i,a_{i+1},a_{i+2},\dots,a_n)
:= (a_1,\dots,a_{i-1},\hskip -10pt \phantom{a}^{a_i} a_{i+1},
a_i,a_{i+2},\dots,a_n)$}$$
(we write $\hskip -6pt \phantom{b}^{a} b$ for the conjugate $aba^{-1}$)
extends to an action of $B_n$ on $A^n$ (when $G=A$ is the free
group $F_n$, we recover the usual action of $B_n$ on $F_n$).
Clearly, the product map
\begin{eqnarray*}
A^n & \longrightarrow & G \\
(a_1,\dots,a_n) & \longmapsto & \prod_{i=1}^n a_i
\end{eqnarray*}
is invariant with respect to this action. In particular,
for any $g\in G$, one has an action of $B_{l_A(g)}$ on
$\Red_A(g)$.

\subsection{A dual Matsumoto property}

In classical Coxeter theory, the Matsumoto property expresses that two
reduced $S$-decompositions of the same $w\in W$
can be transformed one into the other by
successive uses of braid relations. The ``dual braid
relations'' will be defined in the next section -- the next
proposition will then appear as a dual analog of the Matsumoto property.

\begin{prop}
\label{length2}
Let $(W,T)$ be a reflection group.
Let $w\in W$. If $w$ is a parabolic Coxeter element of $(W,T)$,
then the action of $B_{l_T(w)}$ on $\Red_T(w)$ is transitive.
\end{prop}

This proposition indicates why Coxeter elements play a special part
in the dual approach. Indeed, the $B_{l_T(w)}$-action on $\Red_T(w)$
may not be transitive for an arbitrary $w$. Take for example
the Coxeter system of type $B_2$: let $W = < s,t | stst=tsts, s^2=t^2=1>$. The
set $T$ consists of $s$, $t$, $u:=tst$ and $v:=sts$, and we have
$\Red_T(stst) = \{ (s,u), (u,s), (t,v), (v,t) \}$. Since $s$ commutes
with $u$ and $t$ commutes with $v$, the action of $B_2$ has two orbits.

\begin{proof} Thanks to Lemma \ref{paracox},
it is enough to deal with the case of Coxeter elements: indeed,
if $w$ is a parabolic Coxeter in $(W,T)$, it is a
Coxeter element in some $(W_I,T_I)$,
and $\Red_T(w) = \Red_{T_I}(w)$.

We prove the proposition, for Coxeter elements,
by induction on the rank $n$ of $(W,T)$.
It is obvious when $n$ is $0$ or $1$.

Assume $n>1$, and assume the proposition is known for
Coxeter elements in parabolic subgroups of rank $n-1$.
Let $c\in W$ be a Coxeter element in $W$.
Our goal is to prove that $\Red_T(c)$ forms a single orbit under
the action of $B_n$. Denote by $\bullet$ the concatenation of
finite sequences; we write $(t)\bullet \Red_T(tc)$ for
$\{ (t) \bullet u | u\in \Red_T(tc)\}$.
We have
$$\Red_T(c) = \bigcup_{t\in T} (t) \bullet \Red_T(tc).$$

For all $t\in T$, $tc$ is a parabolic Coxeter element. The induction
assumption ensures that the action of $B_{n-1}$ is transitive
on $\Red_T(tc)$. In particular, since the $B_{n-1}$-action
on the last $n-1$ terms is a restriction of the $B_n$-action,
for any $u\in \Red_T(tc)$, the $B_n$-orbit of $(t) \bullet u$ contains
the whole $(t) \bullet \Red_T(tc)$. To conclude, it is enough to exhibit
a particular element of $\Red_T(c)$ such that its orbit under the action
of $B_n$ contains at least one element in each of the
$(t) \bullet \Red_T(tc)$. This may be done as follows:

Let $(L,R)$ be a chromatic pair such that $c=c_{L,R}$.
Write $L=\{s_1,\dots,s_k\}$, $R=\{s_{k+1},\dots,s_n\}$.
Clearly, $(s_1,\dots,s_n) \in \Red_T(c)$.
A direct computation shows that, for all $i\in \{1,\dots,n\}$, the
word
$$\sigma_1^{-1} \sigma_2^{-1} \dots \sigma_{i-1}^{-1} (s_1,\dots,s_n)$$
starts with $s_i$.
Another straightforward computation yields the following:
$$\forall (t_1,\dots,t_n)\in \Red_T(c), (\sigma_{n-1}\dots \sigma_1)^n
(t_1,\dots,t_n) = (ct_1c^{-1},\dots,ct_nc^{-1}).$$
As a consequence, for all $i\in \{1,\dots,n\}$ and for all positive
integer $k$,
$$((\sigma_{n-1}\dots \sigma_1)^{nk}\sigma_1^{-1} \sigma_2^{-1} \dots \sigma_{i-1}^{-1}) (s_1,\dots,s_n)$$ is
an element in the $B_n$-orbit of $(s_1,\dots,s_n)$ starting by
$c^ks_ic^{-k}$. Lemma \ref{orbit} ensures that all elements of $T$
are of the form $c^ks_ic^{-k}$. 
\end{proof}

We have the following immediate consequence.

\begin{coro}
Let $(W,T)$ be a reflection group. Let $w\in W$ be a parabolic
Coxeter element. Let $(t_1,\dots,t_k)\in \Red_T(w)$. The subgroup
$<t_1,\dots,t_k>\subset W$ does not depend on the choice of
$(t_1,\dots,t_k)$ in $\Red_T(w)$.
\end{coro}

In the context of the corollary, we denote the subgroup 
$<t_1,\dots,t_k>$ by $W_w$.
Let $T_w:=W_w \cap T$.
The reflection group $(W_w,T_w)$ is a parabolic subgroup
of $(W,T)$. All parabolic subgroups may be obtained this way.

\begin{defi}
Let $(W,T,c)$ be a dual Coxeter system. A {\sf standard parabolic
subgroup} (with respect to $c$) is a parabolic subgroup
of the form $(W_w,T_w)$, where $w\prec_T c$.
\end{defi}

Let $\rho$ be a realization of $(W,T)$ in $\GL(V_{\BR})$.
For each $w\in W$, we set $K_w := \ker(\rho(w) -\Id)$.

The next proposition summarizes the main results from \cite{bradywatt2};
it is a refinement of Lemma \ref{carter}.

\begin{prop}
\label{brady}
Let $w\in W$.
\begin{itemize}
\item[(1)] For all $w'\in W$, we have $w'\prec_T w 
\Leftrightarrow K_{w'} \supset K_w$.
\item[(2)] Let $w',w''\in W$. Assume that both $w'\prec_T w$ and
$w''\prec_T w$, and that 
$K_{w'} = K_{w''}$. Then $w'=w''$.
\end{itemize}
\end{prop}

In particular, the map $w\mapsto K_w$ is a poset isomorphism from 
$$(\{w\in W | w\prec_T c\},\prec_T)$$ (the underlying set is
denoted by $P_c$ in  \ref{ggroup})
to its image in the subspaces poset of $V_{\BR}$.
To each subspace of $V_{\BR}$ corresponds a parabolic subgroup of $(W,T)$.
Subspaces in the image of $w\mapsto K_w$
correspond to standard parabolic subgroups.

\section{The dual braid monoid}

Throughout this section, we work with a fixed reflection group $(W,T)$.
We often use ``light'' notations, not explicitly refering to $(W,T)$,
though of course all
constructions are depending on $(W,T)$.

\subsection{The dual braid relations}

\begin{defi}
Let c a Coxeter element.
We say that two reflections $s,t\in T$ are {\sf non-crossing} (with
respect to $(W,T,c)$), and we write
$s\parallel_c t$, if $st\prec_T c$ or $ts\prec_T c$. 
\end{defi}

For any $s,t\in T$, the property $s\parallel_c t$ is equivalent
to the existence of an element of $\Red_T(c)$ in which both
$s$ and $t$ appear (use the braid group action from the previous
section).
This relation is symmetric but in general not transitive.
Note that the notion really depends on $c$.

Throughout this section, if $A$ is an alphabet, we call
{\sf relation}  an unordered pair of words in $A^*$. We write
it $u=v$, or equivalently $v=u$, where $u$ and $v$ are the two words.
E.g., in the next definitions, the dual braid relations are identities
between length $2$ words on the alphabet $T$.

\begin{defi}
Let $c$ be a Coxeter element.
The {\sf dual braid relations} (with respect to $(W,T,c)$) 
are the relations of the form
$st= \hskip -6pt \phantom{t}^s t s$, 
where $s,t\in T$ are such that $st\prec_T c$.
\end{defi}

A consequence of Proposition
\ref{brady} (2) is that if we have both $st\prec_T c $ and 
$ts\prec_T c$, then $st=ts$ (use that $K_{st}=K_{ts}$ in any
realization). The  
dual braid relations associated with $st$ and $ts$
are then both equal to $st=ts$.
Therefore dual braid relations are indexed by unordered pairs
of non-crossing reflections.

We use the terminology from subsection \ref{ggroup} of the appendix.
The pair $(W,T)$ is a generated group.
As noted before, since $l_T$ is invariant by conjugacy, 
we have $\forall w,w'\in W,
w\prec_T w' \Leftrightarrow w' \succ_T w$, and all elements of $W$,
and in particular Coxeter elements,
are $T$-balanced.

\begin{defi}
Let $c$ be a Coxeter element.
Let $P_c$ be the pre-monoid of divisors of $c$ in the generated
group $(W,T)$ (as in the appendix).
The corresponding monoid
$\BM(P_c)$ is called the
{\sf dual braid monoid} (with respect to $(W,T,c)$).
\end{defi}

The object of this section is the study of
the combinatorics of $\BM(P_c)$.

We start by deducing from the ``dual Matsumoto property'' that the dual 
braid monoid
is presented by the dual braid relations. In other words, only a small
fragment of the pre-monoid structure of $P_c$ 
is needed to recover all relations in $\BM(P_c)$:

\begin{theo}
\label{ruth}
Let $c$ be a Coxeter element.
The embedding $T\hookrightarrow P_c$ between generating
sets induces an isomorphism
$$< T | \text{dual braid relations} >_{\text{\em monoid}} \;  \;
\simeq \; \BM(P_c).$$
\end{theo}

\begin{proof}
The monoid $\BM(P_c)$ is generated by its atoms, which are exactly the
elements of $T$.
A presentation for $\BM(P_c)$, with respect to this generating set,
is obtained by taking all relations of the form $u=v$, where $u$ and
$v$ are reduced $T$-decompositions of the same $w\in W$, with $w\prec_T c$.
Let us prove that such a relation $u=v$ is a consequence of the dual braid
relations with respect to $(W,T,c)$.
By Lemma \ref{paracox}, the corresponding $w$ is a Coxeter element in a
parabolic subgroup $(W_I,T_I)$. Of course, $u$ and $v$ are
reduced $T_I$-decompositions of $w$.
By Proposition \ref{length2}, any two
reduced $T_I$-decompositions of $w$ are in the same orbit for the braid
group action. By the very definition of this braid group action, 
this ensures that the relation $u=v$ is a consequence of the
dual braid relations with respect to $W_I$ and $w$. Since $w\prec_T c$,
these ``parabolic'' dual braid relations constitute a subset of the set
of dual braid relations with respect to $W$ and $c$.
\end{proof}

\subsection{Dual relations and classical relations}

The goal of this subsection is to prove that the group of fractions of the
dual braid monoid is isomorphic to the Artin group associated to $W$.
In terms of presentations, this means that the dual braid relations are,
as group-defining relations, ``equivalent'' to the classical braid relations.

\begin{lemma}
\label{parallel}
Let $c$ be a Coxeter element. 
Let $s,t\in T$, with $s\neq t$. We denote by $m_{s,t}$ the 
order of $st$ in $W$. The following assertions are equivalent:
\begin{itemize}
\item[(i)]  $s\parallel_c t$.
\item[(ii)]
The classical braid relation 
$$\underbrace{sts \dots}_{m_{s,t}\;\text{terms}} =
\underbrace{tst \dots}_{m_{s,t}\;\text{terms}}$$ is satisfied
in $\BM(P_c)$.
\end{itemize}
\end{lemma}

\begin{proof}
The implication $(ii)\Rightarrow (i)$ is obvious, since if $s$ and $t$
are crossing,
no dual braid relation can be applied to 
$\underbrace{sts \dots}_{m_{s,t}\;\text{terms}}$.

Let us now prove $(i) \Rightarrow (ii)$.
We set $m:=m_{s,t}$.
Without loss of generality, we may assume $st\prec_T c$.
Let $s_1:=t$, $s_2:=s$ and, for $n>2$,
$s_{n+1}:= \hskip -6pt \phantom{s}^{s_n} s_{n-1}$.
We have, for all $n$, $s_{n+m}=s_n$, and
$$s_2s_1=s_3s_2=s_4s_3 = \dots = s_{m-1}s_{m-2}= s_m s_{m-1} = s_1s_m$$
is a sequence of dual braid relations.

Assume $m$ is even.
By multiple use of the above relations, we have
\begin{eqnarray*}
\underbrace{s_2s_1s_2s_1\dots s_2s_1s_2s_1}_{m \; \text{terms} } & = &
\underbrace{(s_1s_m)(s_{m-1}s_{m-2})\dots (s_5s_4)(s_3s_2)}_{m 
\; \text{terms}}\\
& = &
\underbrace{s_1(s_2s_1)(s_2\dots s_1)(s_2s_1)s_2}_{m \; \text{terms}}
\end{eqnarray*}

Assume $m$ is odd. We have
\begin{eqnarray*}
\underbrace{s_2s_1s_2s_1\dots s_2s_1s_2s_1s_2}_{m \; \text{terms} } & = &
\underbrace{(s_1s_m)(s_{m-1}s_{m-2})\dots (s_6s_5)(s_4s_3)s_2}_{m 
\; \text{terms}}\\
& = & \underbrace{s_1(s_2s_1)(s_{2}\dots s_1)(s_2s_1)(s_2s_1)}_{m 
\; \text{terms}}\\
\end{eqnarray*}
\end{proof}

From now on, we fix a chromatic pair
$(L,R)$. Let $S:=L\cup R$, let $c:=c_{L,R}$. We write $L=\{s_1,\dots,s_k\}$,
$R=\{s_{k+1},\dots,s_n\}$, such that $c=s_1\dots s_n$.
In $\BB(W,S)$, we consider the corresponding $\BS$, $\bs_i$ and $\bc$.
We set $$\BT:= \{ \bc^k \bs \bc^{-k} | k\in \BZ, \bs \in \BS\}.$$

The next lemma is a ``braid version'' of Lemma \ref{orbit}.

\begin{lemma}
\label{braidorbit}
Let $\bt,\bu\in\BT$, and let $t,u$ be the corresponding elements of
$W$. For all $m\in \BZ$, we have
$$\bu = \bc^m \bt \bc^{-m}\; \text{\em in $\BB(W,S)$} 
\Leftrightarrow u = c^m t c^{-m} \; \text{\em in $W$}.$$

The canonical morphism $\BB(W,S) \twoheadrightarrow W$ restricts
to a bijection $$\BT \stackrel{\sim}{\rightarrow} T.$$
\end{lemma}

\begin{proof}
It is enough to prove the result when $W$ is irreducible.

The implication $\bu = \bc^m \bt \bc^{-m}\; \text{\em in $\BB(W,S)$} 
\Rightarrow u = c^m t c^{-m} \; \text{\em in $W$}$ is obvious.

Let $s,s'\in S$, let $m,l \in \BZ_{\geq 0} $ such
that $c^m s c^{-m} = c^l s' c^{-l}$ in $W$. To obtain the converse
implication, we have to prove that
$\bc^m \bs \bc^{-m} = \bc^l \bs' \bc^{-l}$ in $\BB(W,S)$. It suffices
to deal with the case $l=0$. From now on, we assume $c^m s c^{-m} =  s'$.

By Lemma \ref{orbit}, this implies that $m$ is a multiple of $h/2$
(of course, if $h$ is odd, then $m$ must be a multiple of $h$).
According to \cite{bourbaki},
Ch. V, $\S$ 6, Ex. 2, p 140, we have $\bc^{h/2}= \bw_0$ when $h$ is even,
and $\bc^h = \bw_0^2$ with no restriction on $h$.
In any case, we have $\bc^m = \bw_0^{2m/h}$. The conjugation by $\bw_0$ 
is a diagram automorphism of the classical diagram for $\BB(W,S)$; the
relation $\bc^m \bs \bc^{-m} = \bs'$  follows immediately.

Since the natural map $\BT \rightarrow T$ is $\BZ/h\BZ$-equivariant (for
the respective conjugacy actions by powers of $\bc$ and $c$),
the description
of $T$ given in Lemma \ref{orbit} and the definition of $\BT$ prove that
$\BT \rightarrow T$ is a bijection.
\end{proof}

\begin{lemma}
\label{conjeasy}
Let $t,u\in T$. For all $m\in \BZ$, we have
$$u c^m = c^m t \; \text{\em in $\BM(P_c)$} 
\Leftrightarrow u = c^m t c^{-m} \; \text{\em in $W$}.$$
\end{lemma}

\begin{proof}
An obvious induction reduces the lemma to the case $m=1$.
Let $t\in T$. Let $u:=ctc^{-1}$.
Let $(t_1,\dots,t_n)\in \Red_T(c)$ such that
$t_1= u$ (Lemma \ref{debut}).
The relation $ut_2\dots t_n=t_2\dots t_n t$ is a consequence of the dual
braid relations.
Thus, in $\BM(P_c)$, we have $uc=uut_2\dots t_n= u t_2 \dots t_n t=
t_2 \dots t_n t t = ct$.
\end{proof}

We call {\sf $\bc$-conjugacy relations} the relations of the form
$\bt (\bs_1\dots \bs_n)^m= (\bs_1\dots \bs_n)^m \bs$ (with $\bt\in \BT$,
$\bs\in \BS$ and $m$ a positive integer) which are
true in $\BB(W,S)$.

The group $\BB(W,S)$ has the presentation
$$\BB(W,S) = < \BS | \text{classical braid relations} >_{\text{group}}.$$
Since the $\bc$-conjugacy relations allow the elements
of $\BT$ to be expressed
as conjugates of the elements of $\BS$, a successive use
of Schreier transformations
introducing the redundant generators $\BT - \BS$ leads
to the following presentation:
$$\BB(W,S) \simeq \left< \BT \left| 
\begin{matrix}
\text{classical braid relations on $\BS$} \\
+ \; \text{$\bc$-conjugacy relations}
\end{matrix}
\right. \right>_{\text{group}}.$$

The group $\BG(P_c)$ has the presentation
$$\BG(P_c) = < T | \text{dual braid relations} >_{\text{group}}.$$
According to the Lemma \ref{conjeasy}, the ``$c$-conjugacy relations''
are consequences of the
dual braid
relations. If $s,s'\in S$, we have $s\parallel_c s'$; by Lemma \ref{parallel},
the classical braid relation involving $s$ and $s'$ is true in
$\BG(P_c)$. Adding these two sets of redundant relations, 
we obtain
$$\BG(P_c) = \left< T \left| 
\begin{matrix}
\text{dual braid relations}\\ 
+ \; \text{$c$-conjugacy relations} \\
+ \; \text{classical braid relations on $S$}
\end{matrix}
\right. \right>_{\text{group}}.$$

This proves that the bijection $\BT \stackrel{\sim}{\rightarrow} T$ from Lemma
\ref{braidorbit} extends to a group morphism
$$\BB(W,S) \rightarrow \BG(P_c).$$
The morphism is invertible, thanks to the following fact:

\begin{fact}
\label{factA}
Let $\bt,\bu \in \BT$, and let $t,u$ be the corresponding elements of
$W$. Assume that $tu\prec_T c$. Let $tu=uv$, with $v\in T$, be the corresponding
dual braid relation.
Then $\bt\bu = \bu \bv$ in $\BB(W,S)$ (where $\bv\in \BT$ corresponds to $v$).
\end{fact}

\begin{proof}
We only have a case-by-case proof.
It is enough to deal with the irreducible case.
The exceptional types are dealt with by computer, using the package
{\sf CHEVIE} of {\sf GAP}.
The dihedral case is obvious.
For the infinite families $A$, $B$ and $D$, see section \ref{ABD}.

Note however that the geometric interpretation of the next section allows
a reformulation of this fact which, we hope, could lead to a general
proof.
\end{proof}

This completes the proof of:

\begin{theo}
\label{mainA}
The bijection $T\stackrel{\sim}{\rightarrow} \BT$ extends to a group
isomorphism $$\BG(P_c)\stackrel{\sim}{\rightarrow} \BB(W,S).$$
\end{theo}

We will later see that $\BM(P_c)$ embeds in $\BG(P_c)$, and therefore
that $\BM(P_c)$ is isomorphic to the submonoid of $\BB(W,S)$ generated
by $\BT$.

For all $s,t\in T$ such that $st\prec_T c$, let us denote by 
$t \stackrel{s}{\rightarrow}\hskip -6pt \phantom{t}^s t$
the dual braid relation relation $st= \hskip -6pt \phantom{t}^s t s$.
Viewing each relation $t \stackrel{s}{\rightarrow}\hskip -6pt \phantom{t}^s t$
as a labelled oriented edge connecting $t$ and $\hskip -6pt \phantom{t}^s t$,
and putting together all dual braid relations with respect to $c$,
we obtain a {\sf labelled oriented graph} with $T$ as set of vertices
(the edges are themselves labelled by vertices).
Together, Theorems \ref{ruth} and \ref{mainA} show that $\BB(W,S)$ admits
what is called a {\sf labelled oriented graph presentation}
or {\sf LOG presentation}. These presentations have been studied by
various authors (see for example \cite{howie}), and are related to 
topological properties.
A typical example of such a presentation is the Wirtinger presentation
for the fundamental group of a link complement. The author thanks Ruth Corran
for pointing out this interpretation.

Since the elements of $\BT$ are conjugates of elements in $\BS$, any
presentation of $\BB(W,S)$ with $\BT$ as set of generators
yields, by addition of quadratic relations, a presentation for $W$.
As a corollary of the above theorem, 
we obtain a ``dual Coxeter presentation'' for $W$.

\begin{coro}
The group $W$ has the following
presentation:
$$W\; \simeq \; \; 
< T | \text{dual braid relations relative to $c$} \; + \; \forall t\in T,t^2=1 
>_{\text{\em group}}.$$
\end{coro}

Instead of deducing this corollary from Theorem \ref{mainA}, we could
have given a direct proof without case-by-case, 
using the trivial analog of Fact \ref{factA}
where the conclusion ``$\bt\bu = \bu \bv$ in $\BB(W,S)$'' is replaced
by ``$tu = uv$ in $W$''.

There is a well-known example of a presentation for $W$ involving
all the reflections: the 
{Steinberg presentation} of $W$, where, in addition to quadratic relations,
all true relations of the form
$st=tu$ are taken (not just those corresponding to non-crossing reflections).
But, when removing the quadratic relations, the Steinberg presentation does
not give a presentation of the braid group.

\subsection{The dual monoid is a Garside monoid}
Since $T$ consists of involutions, word reversing provides
a bijection between $\Red_T(c)$ and $\Red_T(c^{-1})$,
and the posets $(P_c,\prec_T)$ and $(P_{c^{-1}},\succ_T)$ are 
isomorphic. Since $T$ is invariant by conjugacy and the Coxeter
elements $c$ and $c^{-1}$ are conjugate, the posets
$(P_{c^{-1}},\succ_T)$  and $(P_c,\succ_T)$ are isomorphic.
Hence  $(P_c,\prec_T) \simeq (P_c,\succ_T)$ (but, in general,
the identity $P_c\rightarrow P_c$ is not an isomorphism).

\begin{fact}
\label{factB}
The poset $(P_c,\prec_T)$ is a lattice.
\end{fact}

\begin{proof}
Here again, we only have a case-by-case proof (the reduction to the
exceptional case is obvious). We hope that
the geometric approach of the next section will eventually provide 
a general proof.

The exceptional types are dealt with by computer, using {\sf GAP}.
The type $I_2(e)$ is trivial (the poset has height $2$, with only
one maximal element).

For type $A$, $B$ and $D$, see section \ref{ABD}.
Note that, using Theorem \ref{BDM}, we only have to check
that pairs of reflections have a right lcm.
\end{proof}

Applying Theorem \ref{BDM}, we obtain the following:

\begin{theo}
\label{mainB}
The dual braid monoid $\BM(P_c)$ is a Garside monoid.
\end{theo}

Being Garside is a very strong property for a monoid.
Some general properties of Garside monoids are listed at
the end of the appendix (see also section \ref{applications}).
These properties may be seen as a justification for the study
of the dual braid monoid.

\begin{coro}
\label{coroB}
The dual braid monoid $\BM(P_c)$ is isomorphic to the submonoid
of $\BB(W,S)$ generated by $\BT$.
\end{coro}

\begin{proof}
Garside monoids satisfy the embedding property: the natural 
monoid morphism
$\BM(P_c) \rightarrow \BG(P_c)$ is injective.
We conclude using Theorem \ref{mainA}.
\end{proof}

\subsection{Automorphisms of the dual braid monoid}

\begin{defi}
We say that a monoid $M$ is {\sf symmetric} if it admits a generating
set $A$ such that the identity map $A\rightarrow A$ extends to an 
anti-automorphism of $M$.
\end{defi}

Clearly, this is equivalent to the existence of a presentation such
that whenever $u=v$ is a relation, $\overline{u}=\overline{v}$ is also a
relation (where $\overline{u}$ and $\overline{v}$ are the reversed words).

For example, the classical braid monoids are symmetric monoids.

Dual braid monoids are (in general) not symmetric. Consider for example
the reflection group of type $A_2$. The set $T$ consists
of three elements $s,t,u$, such that $st=tu=us$. Choose $c=st$ (the other
choice is $c=ts$). Since $T$ is the set of atoms of $\BM(P_c)$, any generating
set for $\BM(P_c)$ must contain $T$. But the reversed defining 
relation $ts=ut$ is not true in $\BM(P_c)$ (if it were true, then
$ut$ would be a minimal common right multiple of $t$ and $u$,
thus equal to their right lcm $tu$ --- but $ut\neq tu$ in $W$, which is 
a quotient of $\BM(P_c)$).

Nonetheless, $\BM(P_c)$ admits anti-automorphisms:

Let $c$ be Coxeter element. Then $c^{-1}$ is a Coxeter element and
the identity map $T\rightarrow T$ induces
an anti-isomorphism $$\psi_c:\BM(P_c) 
\underset{\op}{\stackrel{\sim}{\rightarrow}} \BM(P_{c^{-1}}).$$

For any $w\in W$, the conjugate $\hskip -6pt \phantom{c}^w c$
is a Coxeter element. The bijection
$T  \rightarrow  T, t  \mapsto  \hskip -6pt \phantom{t}^w t$
extends to an isomorphism
$$\phi_{c,w}: \BM(P_c) \stackrel{\sim}{\rightarrow} 
\BM(P_{\hskip -6pt \phantom{c}^w c}).$$

Let $(L,R)$ be a chromatic pair such that $c=c_{L,R}$.
Then $c^{-1}=c_{R,L}$.
We have $s_R c s_R^{-1} = c^{-1}$.
If $(W,T)$ is irreducible, the centralizer of $c$ in $W$ is the cyclic
subgroup generated by $c$, so any $w\in W$ such that 
$w c w^{-1} = c^{-1}$ is of the form $s_Rc^k$.

We set $$\Theta:= \psi_c^{-1} \circ \phi_{c,s_R}$$ and
 $$\Theta':= \psi_c^{-1} \circ \phi_{c,s_Rc}.$$

By looking at the conjugacy action of $\Theta$ and $\Theta'$ on $T$,
we obtain the following result:

\begin{prop}
Assume $(W,T)$ is irreducible with Coxeter number $h$.
The maps $\Theta$ and $\Theta'$ are involutive anti-automorphisms
of $\BM(P_c)$. They satisfy a classical braid relation of length
$h$. This defines an action of the dihedral group $I_2(h)$ on $\BM(P_c)$, 
such that reflections act by anti-automorphisms and rotations by
automorphisms.

If $ZW$ is trivial, this representation of $I_2(h)$ is faithful. Otherwise,
$ZW$ has order $2$, $h$ is even, and the kernel of the representation
is the center of $I_2(h)$.
\end{prop}

It should not be too difficult to answer the following:

\begin{question}
Let $C$ be the conjugacy class of Coxeter elements in $W$.
Is it possible to find a transitive system of isomorphisms between the
$(\BM(P_c))_{c\in C}$?
\end{question}

\section{Local braid monoids}
\label{geometric}

We give in this section a geometric description of
the dual monoid.
The classical monoid has an interpretation in terms of walls and chambers or,
in other words, in terms of the convex geometry of the hyperplane 
arrangement, seen from a real basepoint. We prove that the dual monoid
has an analogous interpretation, except that one has to look to the 
complexified hyperplane arrangement from a $h$-regular eigenvector.
Hence the dual monoid is indeed a new {\em point of view} on braid groups...

The structure of the section is as follows: \ref{geonotations} and
\ref{lr} contain only generalities; the material in \ref{B0} is more
or less standard, we include it to justify certain computations; 
in \ref{local}, we construct for each basepoint
a ``local'' set of generators and a ``local''
submonoid of the braid group; when the basepoint is a regular eigenvector,
the monoid has certain symmetries, as we will see in \ref{localreg}.
A real basepoint yields the classical monoid. The main results of this
section are in \ref{localdual}, where we interpret the dual monoid
as a certain local monoid.

\subsection{Braid groups}
\label{geonotations}
Let $(W,T)$ be an abstract reflection group. For simplicity, we assume
throughout this section that $(W,T)$ is irreducible.
Let $\rho: W \hookrightarrow \GL(V_{\BR})$ be a realization. We identify
$W$ and $\rho(W)$. We denote by $\CA_{\BR}$ the set of reflecting 
hyperplanes. Let $V$ be the complexified representation
$V\otimes_{\BR}\BC$.
As a complex representation, $V$ is irreducible; as a real representation,
$V=(V_{\BR}\otimes 1 ) \oplus (V_{\BR}\otimes i)$, with
$V_{\BR}\otimes 1 \simeq V_{\BR}\otimes i \simeq V_{\BR}$; we denote by
$\Re$ and $\Im$ the two corresponding $W$-equivariant projections
$V\rightarrow V_{\BR}$.

We denote by $\CA$ the set of reflecting hyperplanes in $V$.
More generally, we often use curly letters for subsets of $\CA$ and plain
letters for the corresponding subsets of $T$ (we preferred not to change
the standard notation $\CA$ into $\mathcal{T}$). We also use curly letters
for chambers: a {\sf (real) chamber} $\CC$ is a connected component
of $V_{\BR} - \bigcup_{H_{\BR} \in \CA_{\BR}} H_{\BR}$. To a chamber
$\CC$, we associate the set $\CS\subset \CA$ of walls of $\BC$.
The corresponding $S\subset T$ is such that $(W,S)$ is a Coxeter system.

We set $$V^{\reg} := V -\bigcup_{H\in \CA} H.$$
The covering $V^{\reg} \twoheadrightarrow W\backslash V^{\reg}$ is
unramified (since $W$ acts, as usual, on the left on $V$, the notation
$W \backslash V^{\reg}$ for the orbit space is more acurate than
the usual $V^{\reg}/W$).

The {\sf braid group } of $(W,T)$ is the fundamental
group of $W\backslash V^{\reg}$. Of course, this is well-defined only
up to the choice of a basepoint -- this choice will appear to be crucial 
here.

{\flushleft \bf Remark on the terminology.}
The term ``Artin group'' should only be used to refer to the abstractly
presented group $\BB(W,S)$. Of course, by a result of Brieskorn
(see Theorem \ref{brieskorn} below),
$\BB(W,S)$ is isomorphic to the braid group of $W$. But the isomorphism
is not canonical.

{\small
\subsection{Left and right, morphisms and anti-morphisms}
\label{lr}

Suppose given a Galois covering $p:X \twoheadrightarrow G\backslash X$.
Let $x\in X$ and let $\gamma$ be a loop in the pointed space 
$(G \backslash X,p(x))$.
The homotopy lifting of $\gamma$ is a path $\tilde{\gamma}$ in $X$ from
$x$ to some $y\in Gx$. Since $p$ is Galois, there
is a unique $g\in G$ such that $y=gx$.
This defines a map $\pi_1(G\backslash X,p(x)) \rightarrow G$.
The standard rule for concatenation of paths ($\gamma\gamma'$ is defined
when the ending point of $\gamma$ is the starting point of $\gamma'$)
makes this map an anti-morphism, rather than a morphism.

Both the classical monoid and the dual monoid are endowed with canonical
morphisms to $W$. Hence, if we want to compare them with the braid
group, we should look for submonoids of the braid group which are
anti-isomorphic to them.
This is usually not an issue with the classical braid monoid since it is
a symmetric monoid. But the dual monoid is not symmetric. 
It will be
important here to keep in mind that the natural formulation for
Brieskorn's theorem (\ref{brieskorn}) involves an anti-isomorphism.

{\bf \flushleft Remark.}
To get rid of this ``anti'', it would have been possible to modify
one of our conventions:
\begin{itemize}
\item Replace the above map $\gamma \mapsto g$ by $\gamma \mapsto g^{-1}$.
The solution only pushes the ``anti'' further:
it makes the canonical map from the classical (or dual) monoid to $W$
an anti-morphism.
\item Invert the concatenation rule for paths, so that the
fundamental groupoid of $X$ acts on the left on $X$. This solves all the 
problems (but violates a quasi-universal convention).
\item Make $G$ act on the right on $X$.
\end{itemize} }

\subsection{The Brieskorn basepoint}
\label{B0}
We describe here a couple of tricks for computing in braid groups,
inspired by \cite{brieskorn}. In type $A$, given a loop representing
a braid, one may write down a word by looking at where and how the
strings cross in the real projection. This can be generalized to all
types.

Let $\CC$ be a real chamber, with set of walls $\CS$.
We choose a left/right colouring of the Coxeter diagram on $S$. We have a
corresponding colouring of $\CS$. We denote by $\CL$ the set of 
left walls, by $\CR$ the set of right walls.

Since we want our discussion to include types $H$ and $I$, we use a
standard substitute for root systems.
For each $H\in \CA$, we fix a linear form $l_H: V_\BR \rightarrow \BR$
with kernel $H_\BR$, and such that $\forall x\in \CC, l_H(x) > 0$.
By extension of scalars, we view $l_H$ as a linear form $V \rightarrow \BC$
with kernel $H$. It is uniquely defined, up to multiplication by an
element of $\BR_+^*$ (we could normalize $l_H$ using the invariant 
scalar product on $V_{\BR}$, but this is not crucial here).
The following conditions
are equivalent, for a given $v\in V_\BR$:
\begin{itemize}
\item[(i)] The vector $v$ is in $\CC$.
\item[(ii)] For all $H\in \CA$, we have $l_H(v)>0$.
\item[(iii)] For all $H\in \CS$, we have $l_H(v)>0$.
\end{itemize}
Since $(l_H)_{H\in\CS}$ is a basis of the dual of $V_\BR$, for any $H'\in\CA$,
we have $l_H'=\sum_{H\in\CS} \alpha_{H',H} l_H$. A consequence of
$(ii) \Rightarrow (iii)$ is that all coefficients are in $\BR_{\geq 0}$.

The space $\Re^{-1}(\CC) = \CC\otimes 1 + V_\BR\otimes i \subset V^{\reg}$ is
contractible. Thus we may choose it as a basepoint for $V^{\reg}$.
More precisely, for any $v\in \Re^{-1}(\CC)$, the homotopy exact
sequence of the triad $\{v\} \subset \Re^{-1}(\CC) \subset V^{\reg}$ yields a 
canonical isomorphism 
$$\pi_1(V^{\reg},v) \simeq \pi_1(V^{\reg},\Re^{-1}(\CC)).$$
Denote by $p$ the quotient map $V^{\reg} \rightarrow W \backslash V^{\reg}$.
The space $B_0:=p(\Re^{-1}(\CC))$ is contractible and
can be used as a ``basepoint'' for $W \backslash V^{\reg}$.
We call $B_0$ the {\em Brieskorn basepoint} of $W \backslash V^{\reg}$.
For any $w,w'\in W$ and  any
path $\gamma$ in $V^{\reg}$ such
that $\Re(\gamma(0))\in w\CC$ and $\Re(\gamma(1))\in w'\CC$
unambiguously defines an element of
$\pi_1(W \backslash V^{\reg},B_0)$,
the latter group being, for any $x_0\in B_0$,
canonically isomorphic to $\pi_1(W \backslash V^{\reg},x_0)$.

\begin{defi}
Let $\gamma$ be a differentiable path in $[a,b]\rightarrow V^{\reg}$.
We say that $t\in [a,b]$ is a {\sf critical time} for $\gamma$ if
$p(\gamma(t)) \notin B_0$.
Let $C_{\gamma}$ be the set of critical times.
We say that $\gamma$ is {\sf non-singular} if all three conditions hold:
\begin{itemize}
\item[(a)] We have $a\notin C_\gamma$ and $b\notin C_\gamma$.
\item[(b)] The set $C_\gamma$ is finite.
\item[(c)] For each $t\in C_\gamma$, there is a unique $H_t\in \CA$ such that
$l_{H_t}(\Re(\gamma(t))) = 0$, and the tangent line to $\Re\circ\gamma$ at $t$
is not included in $\Re(H_t)$.
\end{itemize}
\end{defi}

Condition (b) actually follows from (c), which could be rephrased as
``$\Re\circ\gamma$ is transverse to each stratum of the real hyperplane
arrangement''.

This notion allows a practical reformulation of the
main results in \cite{brieskorn}. Though not explicitly stated by
Brieskorn, this reformulation follows easily from his construction.
We leave the details to the reader.

\begin{theo}[after Brieskorn]
\label{brieskorn}
There exists a (unique, generating)
subset $(\bs_H)_{H\in \CS}$ of $\pi_1(W \backslash V^{\reg},B_0)$
such that, for
 any non-singular differentiable path
 $\gamma:[0,1]\rightarrow V^{\reg}$ such that
 \begin{itemize}
\item $\gamma(0)\in \Re^{-1}(\CC)$,
\item $\gamma$ has a unique critical time $t_0$,
\end{itemize}
if we denote by $H_0$ the hyperplane such that $l_{H_0}(\Re(\gamma(t_0)))=0$
($H_0$ is always a wall of $\CC$), we have:
\begin{itemize}
\item if $\im(l_H(\gamma(t_0))) > 0$, then $\gamma$ represents $\bs_{H_0}$,
\item if $\im(l_H(\gamma(t_0))) < 0$, then $\gamma$ represents $\bs_{H_0}^{-1}$.
\end{itemize}
These generators realize an explicit anti-isomorphism
$$\pi_1(W \backslash V^{\reg},B_0) \underset{\op}{\simeq} \BB(W,S).$$
\end{theo}

The practical aspect of this reformulation is in the
following explicit recipe
for translating non-singular paths
into elements of the Artin group. Let $\gamma:[a,b] \rightarrow V^{\reg}$
be a non-singular differentiable path.
Start by ordering $$t_1<t_2<\dots < t_k$$ the critical times.
Let $a_0=a < a_1 < \dots < a_{k-1} < a_k=b$ be such that
$$a_0<t_1<a_1<t_2<\dots < a_{k-1} < t_k < a_k.$$
For $i=1\dots k$, we denote by $\gamma_i$ the restriction
of $\gamma$ to $[a_{i-1},a_i]$. 
We have $\gamma = \gamma_k\gamma_{k-1}\dots \gamma_2\gamma_1$ 
(remember we have an unusual
convention for concatenating paths).
Hence we are left with the problem of determining the image
of a given $\gamma_i$.
We have $\gamma_i(a_{i-1}) \in p^{-1}(B_0) = \bigcup_{w\in W} \Re^{-1}(w\CC)$.
Let $w_i$ be the unique element of $W$ such that $\gamma_i(a_{i-1})\in
\Re^{-1}(w_i\CC)$. The path $w_i^{-1}\gamma_i$ represents the same
element of $\pi_1(W \backslash V^{\reg},B_0)$ as $\gamma$. Since the 
initial point of $w_i^{-1}\gamma_i$  is in $\Re^{-1}(\CC)$, it is determined
according to Theorem \ref{brieskorn}.

What can we do with a singular path $\gamma: [0,1]\rightarrow V^{\reg}$?
If the endpoints of $\gamma$
are not in $B_0$, then the real projection is {\bf really} ambiguous, since
$\gamma$ is not a relative loop in the ``pointed'' space $(W \backslash V^{\reg},B_0)$.
If the endpoints are indeed in $B_0$, then 
we may always find $\gamma'$ non-singular in the homotopy class of
$\gamma$. Being non-singular is actually a ``generic'' property,
in the sense that one may desingularize $\gamma$ by arbitrary small
perturbations.

(Alternatively, 
desingularization could be avoided by replacing the rudimentary recipe
given above
by a more sophisticated one, able to handle certain paths crossing more than
one real hyperplane at a time.)

\subsection{Local monoids}
\label{local}
For $v,v'\in V$, we denote by $[v,v']$ the affine segment between $v$ and $v'$
(in other words, the convex hull of $\{v,v'\}$).

\begin{defi}
Let $v\in V^{\reg}$. We say that an hyperplane $H\in \CA$ is {\sf visible
from $v$} if and only $$\forall H'\in \CA, [v,s_Hv]\cap H' \neq \varnothing
\Rightarrow H=H'.$$

We set $\Vis_v:=\{ H\in \CA | H \; \text{is visible from $v$}\}$.
\end{defi}

Assume we are given, for each $H\in\CA$, a linear form $l_H$ with
kernel $H$. Clearly: ``$H$ is visible from $v$'' $\Leftrightarrow$
``$\forall H'\in\CA-\{H\}, 0 \notin [l_{H'}(v),l_{H'}(s_Hv)]$''.

Let $v\in V^{\reg}$, with image $x$ in $W \backslash V^{\reg}$.
Let $H\in\CA$. Assume that $H$ is visible from $v$.
Then the path
$$\gamma: [0,1] \rightarrow V^{\reg}, 
t \mapsto (1-t) v + t \frac{v+s_H(v)}{2}$$
is a path from $v$ to $H$ in $V^{\reg}$ (in the sense of \cite{zariski},
section 2.1). The composition $\overline{\gamma}$
of $\gamma$ with the quotient map
$V^{\reg} \twoheadrightarrow W \backslash V^{\reg}$ is a path
from $\overline{v}$ to the discriminant.
As explained in \cite{zariski}, section 2.1, this path defines a
generator-of-the-monodromy in $\pi_1(W \backslash V^{\reg},x)$.
Let us denote by $\bs_{v,H}$ this generator-of-the-monodromy.

\begin{defi}
Let $v\in V^{\reg}$. The elements of $\{ \bs_{v,H} | H \in \Vis_v\}$
are called {\sf local generators} at $v$.
The {\sf local (braid) monoid at $v$}, denoted by $M_v$,
is the submonoid of $\pi_1(W \backslash V^{\reg},p(v))$ generated by
all the local generators.
\end{defi}

For any $v\in V^{\reg}$ and any $w\in W$, we clearly have
$\Vis_{wv} = w\Vis_v$, and $\forall H\in \Vis_v, \bs_{v,H}=\bs_{wv,wH}$.
Thus $M_v = M_{wv}$.
Let $x:=p(v)$. We set
$M_x:=M_v$ call it the {\sf local monoid at $x$}. It
does not depend on the choice of $v$ in $p^{-1}(x)$.

{\bf \flushleft Remark.}
The different
visibility conditions define a certain stratification of $V^{\reg}$.
The maximal strata, with real dimension $\dim_{\BR}(V)$, are those
on which all hyperplanes are visible. The structure of the local monoid
only depends on the stratum of $v$.
Generically, a point $v\in V^{\reg}$ is in a maximal stratum,
all hyperplanes
are visible from $v$ and the structure of the local
monoid is stable by small modification of $v$. When some hyperplanes
are invisible from $v$, the local monoid is not stable.
Note also that the structure of the local monoid
does not only depend on which hyperplanes are visible: ``how'' they are
visible is important. 
For example, the dual monoid will appear later to be
a particular example of local monoid, corresponding a certain
maximum stratum, but there are maximal strata
such that the corresponding local monoid is not the dual monoid.

\smallskip

We start by rephrasing Brieskorn's theorem in terms of local monoids.
When the basepoint is chosen in a real chamber, the local monoid is the
classical positive braid monoid:

\begin{prop}
\label{localclassical}
Let $v\in V^{\reg}$. Let $\CC$ be a real chamber, with set of walls $\CS$.
Assume that $\Re(v)\in \CC$; using Brieskorn's basepoint, we anti-identify
$\pi_1(W\backslash V^{\reg},p(v))$ with $\BB(W,S)$.
\begin{itemize}
\item[(1)] Any $H\in \CS$ is visible
from $v$, and $\bs_{v,H} = \bs_H$.
\item[(2)] Assume that $\Im(v)=0$. Then $\Vis_v = \CS$,
and the anti-isomorphism 
$\pi_1(W\backslash V^{\reg},p(v)) \underset{\op}{\simeq}\BB(W,S)$
restricts to an anti-isomorphism
$$M_v \underset{\op}{\simeq} \BB_+(W,S).$$
\end{itemize}
\end{prop}

\begin{proof}
(1). Let $H\in \CS$. 
Since all hyperplanes have real equations, we have for all $H'\in \CA$
$$[v,s_Hv] \cap H' \neq \varnothing \Rightarrow 
[\Re(v),\Re(s_Hv)] \cap H'_{\BR} \neq \varnothing.$$
The chamber $s_H\CC$ is separated from $\CC$ by
only one wall, $H$. Thus the segment $[\Re(v),\Re(s_Hv)]$ intersects
only one real hyperplane, $H_{\BR}$. 
This proves that $H \in \Vis_v$.
The identity $\bs_{v,H} = \bs_H$ is easy: choose a path representing
$\bs_{v,H}$ and use Theorem \ref{brieskorn}.

(2). If $\Im(v)=0$, then for any $H,H'\in \CA$, we have
$$[v,s_Hv]\cap H' = \varnothing \Leftrightarrow 
[\Re(v),\Re(s_Hv)] \cap H'_{\BR} = \varnothing.$$
If $H\notin \CS$, then the path $[\Re(v),\Re(s_Hv)]$, which exits the 
chamber $\CC$, must cross at least a wall of $\CC$, and $H\notin \Vis_v$.
The second part of the statement follows immediately.
\end{proof}

\begin{lemma}
\label{lambda}
Let $L$ be a complex line (through the origin) in $V$.
Let $L^*:=L-\{0\}$. Assume $L^* \subset V^{\reg}$.
\begin{itemize}
\item[(1)] Let $v,v' \in L^*$. Let $\gamma$ be a path in $L^*$
starting at $v$ and ending at $v'$. The corresponding
isomorphism $$\phi_{\gamma}:
\pi_1(W \backslash V^{\reg},p(v))\stackrel{\sim}{\longrightarrow}
\pi_1(W \backslash V^{\reg},p(v'))$$ does not depend on the choice of $\gamma$.
Let us denote it by $\phi_{v,v'}$. The family
$(\phi_{v,v'})_{v,v'\in L^*}$ is a transitive system
of isomorphisms between the $(\pi_1(W \backslash V^{\reg},p(v)))_{v\in L^*}$.
\item[(2)]
Let $v,v' \in L^*$. Let $\phi_{v,v'}$ be the corresponding isomorphism, 
as in (1).
We have $\Vis_{v'}=\Vis_v$, and
$$\forall H\in \Vis_v, \phi_{v,v'}(\bs_{v,H}) = \bs_{v',H}.$$
The family $(\phi_{v,v'})_{v,v'\in L^*}$ induces by restriction
a transitive system
of isomorphisms between the $(M_v)_{v\in L^*}$.
\end{itemize}
\end{lemma}

The concrete meaning of the lemma is that it makes sense to use
the notations
$\bs_{L^*,H}$ and $M_{L^*}$.

\begin{proof}
(1). {\em A priori}, 
the isomorphism $\phi_\gamma$ only depends on the homotopy class of $\gamma$.
To prove that it does not depend on $\gamma$, it is enough to check
it when $v=v'$, {\ie}, to prove that the conjugacy action
of $\pi_1(L^*, v)$ on $\pi_1(W \backslash V^{\reg}, p(v))$ is trivial.
But $\pi_1(L^*, v)$ is cyclic, generated by an element which is well-known
to be central in $\pi_1(W \backslash V^{\reg}, p(v))$ (see for example
\cite{brmarou}, Lemma 2.4).
The transitivity of the system of isomorphisms follows from the independence
of the choice of $\gamma$.

(2). The visibility condition is invariant by scalar multiplication.
The rest is an easy computation.
\end{proof}

\begin{prop}
\label{localgen}
Let $v\in V^{\reg}$. The group $\pi_1(W \backslash V^{\reg},p(v))$ is generated
(as a group) by the local generators at $v$. 
\end{prop}

\begin{proof}
By Lemma \ref{lambda}, if $\lambda\in\BC^*$, then 
``$\pi_1(W \backslash V^{\reg},p(v))$ is generated (as a group) by the
local generators at $v$''  is equivalent 
``$\pi_1(W \backslash V^{\reg},p(\lambda v))$ is generated (as a group) by the
local generators at $\lambda v$''.

Since $\CA$ is finite, it is always possible to find $\lambda\in\BC^*$
such that $\forall H\in \CA, \re (l_H(\lambda v)) \neq 0$ or, in other
words, $\Re(\lambda v)$ is in a chamber $\CC$. 

By Proposition \ref{localclassical} (1), the set of local generators
at $\lambda v$ contain a classical Artin-type generating subset for
$\pi_1(W \backslash V^{\reg},B_0)$.
\end{proof}

\subsection{Local monoids and regular elements}
\label{localreg}

A regular element in $W$ is an element which has an eigenvector in $V^{\reg}$.
The connection between regular elements and finite order
automorphisms of braid groups was first noticed in \cite{brmi}.

\begin{prop}
\label{localregprop}
Let $w$ be a regular element of $W$, of order $d$.
If the center of $W$ is non-trivial and $d$ is even, 
set $d':=d/2$; otherwise set $d':=d$.
Let $v$ be a regular eigenvector for $w$.

The set $\Vis_v$ is stable by the action of $w$, and
the local monoid $M_v$ admits an automorphism $\phi$ of order $d'$,
such that $$\forall H\in \Vis_v, \phi(\bs_{v,H}) = \bs_{v,wH}.$$
\end{prop}

\begin{proof}
By assumption, we have $wv=\zeta v$, where $\zeta$ is a primitive $d$-th
root of unity. Write $\zeta=e^{2i\pi k/d}$.
Applying Lemma \ref{lambda} to the path 
$\gamma:[0,1] \rightarrow e^{-2i\pi tk/d}$, we obtain an isomorphism
$$\phi: \pi_1(W \backslash V^{\reg},p(v))  \stackrel{\sim}{\longrightarrow}  
\pi_1(W \backslash V^{\reg},p(\zeta^{-1}v)) = \pi_1(W \backslash V^{\reg},p(v))$$
such that, whenever $\forall H\in \Vis_v$,
$$\bs_{v,H} \mapsto \bs_{\zeta^{-1} v,H}=\bs_{w^{-1}v,H}=\bs_{v,wH}.$$
In particular, $\phi$ restricts to an automorphism of $M_v$.

The order of $\phi$ is the same as the order of the action $w$ on
$\Vis_v$; this action is isomorphic to the conjugation action of $w$ on
$S_v:=\{ s_H | H \in \Vis_v\}$. By Proposition \ref{localgen}, the set $S_v$ 
generates $W$. Thus the order of $\phi$ is the smallest $k>1$ such
that $w^k$ is central in $W$. If $ZW=1$, then $k=d$. Otherwise,
the only non-trivial central element is the (unique)
regular element of order $2$.
The conclusion follows.
\end{proof}

\subsection{The dual monoid as a local monoid}
\label{localdual}
This subsection is devoted to the proof of the 
following theorem, which is an analog of Proposition \ref{localclassical} (2)
for the dual monoid.

\begin{theo}
\label{conjecture}
Let $\CC$ be a chamber of the real arrangement with set of
walls $\CS$. Decompose the corresponding $S$ in a chromatic pair
$L\cup R$; we have the corresponding partition $\CS =\CL \cup \CR$.
Let $v$ be a non-zero $e^{2i\pi/h}$-eigenvector for $c:=c_{L,R}$.
Then all hyperplanes are visible from $v$, and the 
assignment
$$ \forall H \in \CA, \bs_{v,H} \longmapsto s_H$$
extends to a unique monoid anti-isomorphism
$$\xymatrix@1{ M_{v} 
\ar[r]^{\sim\phantom{mn}}_{\op\phantom{mn}}  & \BM(P_{c}).}$$
\end{theo}

{\bf \flushleft Remark.} The space $\ker(c-e^{2i\pi/h} \Id)$ is a complex
line (since $a(h)=1$, in the notations of \cite{springer} 3.4 (i)).
The different spaces $\ker(c-e^{2i\pi/h} \Id)$ 
corresponding to different choices of $c$,
are transitively permuted by the action of $W$ (see \cite{springer}, 3.4 (iii)).
These observations immediately imply that the structure of $M_x$ does not 
depend on the choice of $c$ and $v\in \ker(c-e^{2i\pi/h} \Id) \cap V^{\reg}$. 

We fix
$\CC$, $\CS$, $\CL$, $\CR$, $S$, $L$ and $R$ as in the theorem.

The next proposition is a refinement, for Coxeter elements, of a general
remark by Springer (\cite{springer}, bottom of page 173).
We use the notation $\arg$ for the standard
retraction from $\BC^*$ to the unit circle $S^1$.

\begin{prop}
\label{2hgon}
Let $v\in \ker(c_{L,R}-e^{2i\pi / h}) \cap V^{\reg}$.
Consider the map
\begin{eqnarray*}
\theta:\CA & \longrightarrow & S^1 \\
H & \longmapsto & \arg(l_H(v))
\end{eqnarray*}
\begin{itemize}
\item[(1)] The partition $\CS = \CL \cup \CR$  can
be recovered from $\theta$, in the following way:
when $H'\in \CL$ and $H''\in \CR$,
we have
$$\theta(H')/\theta(H'') = e^{i\frac{h-1}{h} \pi}.$$
In particular, $\theta(\CC)$ consists of
exactly two points, at angle $\frac{h-1}{h} \pi$.
\item[(2)] The image $\theta(\CA)$ consists of $h$ consecutive points
on a regular $2h$-gon.
\end{itemize}
\end{prop}

\begin{coro}
The intersection $\ker(c_{L,R}-e^{2i\pi / h}) \cap \Re^{-1}(\CC)$
is non-empty.
\end{coro}

\begin{proof}
Since $\ker(c_{L,R}-e^{2i\pi / h})$ has complex dimension $1$, and since 
the claimed properties are invariant under multiplication
of $v$ by a non-zero complex number, we only have to prove the 
proposition for a particular $v$. It is easy to build one
from the information provided by Bourbaki.

Let us summarize
various results from pages 118--120 in \cite{bourbaki},  Ch. V, \parag 6. 
According to Bourbaki,
it is possible to find $z',z''\in V_\BR$ such that:
\begin{itemize}
\item For any left wall $H'$ and any right wall $H''$, we have
$$\phantom{mmm} l_{H'}(z')=0, \quad l_{H''}(z')>0, \quad l_{H'}(z'')>0 \; 
\; \text{and} \; \; l_{H''}(z'')=0.$$ 
\item The $\BR$-plane $P$ generated by $z'$ and $z''$ is stable by
$s_L$ and $s_R$.
\item The element $s_L$ (resp. $s_R$) acts on $P$ as a reflection
with hyperplane $\BR z'$ (resp. $\BR z''$). Note that there is a unique
(up to scalar multiplication)
scalar product on $P$ invariant by $s_L$ and $s_R$ 
and therefore there is a well-defined notion of angle in $P$.
We have $(\widehat{z'',z'}) = \pi/h$.

\end{itemize}

Since the conditions specifying $z'$ and $z''$ are stable by multiplication
by an element of $\BR_+^*$, we may assume that both of their norms are $1$.
The vector $$n:=\frac{z' - z'' \cos \pi /h }{\sin \pi /h }$$ is
such that $(z'',n)$ is an orthonormal basis.
A direct computation shows that the element $v\in V$ defined by
$$v:= z''\otimes 1 - n \otimes i$$
is an $e^{2i\pi/h}$--eigenvector for $c_{L,R}=s_Ls_R$.

Assume $H'$ is a left wall and $H''$ is a right wall. We have
\begin{eqnarray*}
\frac{\arg(l_{H'}(v))}{\arg(l_{H''}(v))} = \arg\left(
\frac{l_{H'}(v)}{l_{H''}(v)} \right) &
=  &\arg\left(
\frac{l_{H'}(z'') - i l_{H'}(n)}{l_{H''}(z'')- i l_{H''} (n)} \right) \\
& = & \arg\left( 
 \frac{l_{H'}(z'')+il_{H'}(z'')\cot \pi/h}{-i l_{H''}(z')/\sin \pi/n} \right)\\
& = & \arg\left( 
\frac{l_{H'}(z'')}{l_{H''}(z')}\frac{\sin\pi/h +i\cos\pi/h}{-i} \right) \\
& = & \arg\left( \frac{\sin\pi/h +i\cos\pi/h}{-i} \right) \\
& = & -\cos\pi/h + i \sin\pi/h \\
& = & e^{i\pi \frac{h-1}{h}}.
\end{eqnarray*}
This proves claim (1).

(2).
Let $H\in \CA$. For any $w\in W$, the linear form
$wl_H:x\mapsto l_H(w^{-1}x)$ has kernel $wH$, thus
$wl_H = \lambda l_{wH}$, with $\lambda\in\BC$.
In particular, for any $k\in\BZ$, we have
\begin{eqnarray*}
\theta(c_{L,R}^k H)=\arg(l_{c_{L,R}^k H}(v)) & = & 
\arg(c_{L,R}^k l_H( v )) \\ 
& = & 
\pm \arg( l_H (c_{L,R}^{-k} v))\\
& = & \pm\arg(l_H( e^{-2i\pi \frac{k}{h}} v))\\
& = & \pm e^{-2i\pi \frac{k}{h}} \theta(H)
\end{eqnarray*}
By Lemma \ref{orbit}, $T$ is the closure of $S$ for the
conjugacy action of $c_{L,R}$; rephrased in terms of hyperplanes,
this says that $\CA$ is the closure of $\CS$ for the multiplication 
action of $c_{L,R}$. 
Using (1), we see that $\theta(\CA)\cup (-\theta(\CA))$ is the regular
$2h$-gon containing $\theta(\CS)$. Since all $l_H$ are linear
combinations of the $(l_{H'})_{H'\in \CS}$ with real positive coefficients,
$\theta(\CA)$ must consist of the $h$ consecutive points from
$\theta(\CR)$ to $\theta(\CL)$.
\end{proof}

To simplify notations, we now work with an eigenvector
$$v \in \ker(c_{L,R}-e^{2i\pi/d}\Id) \cap \Re^{-1}(\CC)$$ such that,
when $H'\in \CL$ and $H''\in \CR$, one has

$$\arg(l_{H'}(v)) = e^{i\pi \frac{h-1}{2h}} \quad \text{and} \quad
\arg(l_{H''}(v)) =  e^{-i\pi \frac{h-1}{2h}}$$
(the existence of such a $v$ is a consequence of the previous proposition).

We anti-identify $\pi_1(W \backslash V^{\reg},B_0)$ with 
the Artin group $\BB(W,S)$
via Brieskorn theorem.
By Proposition \ref{localclassical} (1), 
$\CS \in \Vis_v$, and for all $H\in \CS$, $\bs_{v,H}= \bs_H$.

\begin{lemma}
\begin{itemize}
\item[(i)]The element $\bc=\bc_{L,R}=
\bs_L \bs_R \in \pi_1(W \backslash V^{\reg},B_0)$ 
is represented by the path $\gamma:[0,1]\rightarrow V^{\reg}, t\mapsto
ve^{2i\pi t/h}$.
\item[(ii)] For all $H\in\CA$ and all $k\in \BZ$,
we have $$\bc^k \bs_{v,H} \bc^{-k} = \bs_{v,c^k H}.$$
\end{itemize}
\end{lemma}

\begin{proof}
(i) follows from Proposition \ref{2hgon} by an immediate computation,
done with the technique described in Subsection \ref{B0}; we leave the 
details to the reader.

From (i) and the proof of Proposition \ref{localregprop}, it follows
that the automorphism $\phi$ from Proposition \ref{localregprop} is
the conjugation by $\bc$. Assertion (ii) follows.
\end{proof}

An immediate consequence of the lemma is that we have
a geometric interpretation of the set $\BT$ defined in the previous
section:

\begin{prop}
Via the anti-identification $$\pi_1(W \backslash V^{\reg},p(v))
\underset{\op}{\simeq} \BB(W,S),$$ the set of
local generators at $v$ coincides with $\BT$.
\end{prop}

Note that we did not use Fact \ref{factA}, nor any case-by-case 
argument, to prove the last proposition. The proposition provides
a geometric setting to check Fact \ref{factA}. In the next section,
we indicate how to do it for types $A$, $B$ and
$D$; the dual braid relations between the elements of $\BT$
will appear to be particular
Sergiescu relations (or, 
for type $D$, some analogs of Sergiescu relations).

Theorem \ref{conjecture} follows from the last proposition and
Corollary \ref{coroB}.

\subsection{Are there other Garside local monoids?}
Ko and Han (\cite{kohan}) have studied a certain class of submonoids of
the type $A$ braid groups. As this class contains all local monoids, their
main theorem has the following consequence
($X_n$ denotes the set of subsets of $\BC$ of cardinal $n$, which is canonically
homeomorphic the space $W \backslash V^{\reg}$, where $W$ is the Coxeter group
if type $A_{n-1}$):

\begin{theo}[after Ko-Han] Let $n\in \BZ_{\geq 1}$ and $x\in X_n$.
If the local monoid $M_x$ is a Garside monoid,
then 
$x$ is included in an affine line or is the set of vertices of a strictly
convex polygon.
\end{theo}

In other words, in the type $A$ case, the classical monoid and the dual 
monoid are the only local monoids which are Garside monoids.

\section{The dual geometries of types $A$, $B$ and $D$}
\label{ABD}

The previous section provides a geometric framework to study
the dual monoid. When $W$ is of type $A$, $B$ and $D$, this
framework can be used to prove Facts \ref{mainA} and \ref{mainB}.

\subsection{Type $A$}
The type $A$ dual monoid coincides with the Birman-Ko-Lee monoid (\cite{BKL}).
In \cite{BDM}, we gave a geometric interpretation of this monoid, 
via what we called ``non-obstructing'' partitions (an equivalent
approach is given in \cite{brady} -- the ``non-obstructing''
partitions from \cite{BDM} are of course as the ``non-crossing''
partitions from \cite{brady} and \cite{reiner}); this interpretation
can be seen as a particular case of the general one given in
section \ref{geometric}. 

Instead of just quoting \cite{BDM} for the lattice property
(Fact \ref{mainB}), we give a survey of the main results, since
they provide an intuitive illustration of some results from the
previous sections.
Formal definitions and
complete proofs can be found in \cite{BDM}.

Let $W$ be the symmetric group ${\mathfrak{S}}_n$, let $T\subset W$ be the
subset of all transpositions. The group $(W,T)$ is an abstract
reflection group of type $A_{n-1}$. The Coxeter elements are the
$n$-cycles.
We choose the standard monomial realization.
The space $W \backslash V^{\reg}$ (see section \ref{geometric}) can be identified
with the space $X_n$ of subsets of $\BC$ of cardinal $n$.
The fiber of $V^{\reg}$ above $x\in X$ is indexed by the total orderings
of $x$: a $n$-tuple $(x_1,\dots,x_n)\in\BC^n$ is above $x$ if
and only if $\{x_1,\dots,x_n\} = x$. 

Let $\mu_n\in X_n$ be the set of complex $n$-th roots of unity.
To fix notations, we choose $v=(e^{2i\pi\frac{1}{n}},e^{2i\pi\frac{2}{n}},
\dots,e^{2i\pi\frac{n}{n}})$. It lies in the fiber over $\mu_n$.
The vector $v$ is a regular
$e^{2i\pi/n}$-eigenvector for the $n$-cycle $c:=(1\; 2 \; \dots \; n)$.
To a transposition $(i\; j)$, we associate the braid $\bs_{i,j}$
represented as follows, by a path where only the $i$-th
and $j$-th strings are moving:
$$
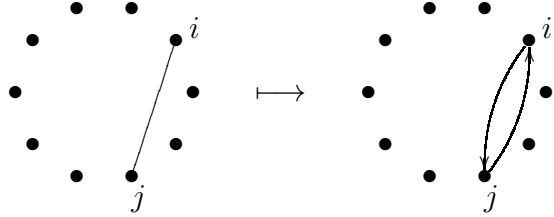
\xy/r2.8pc/:,{\xypolygon10"A"{~>{}}},
{\xypolygon10"B"{~:{(1.25,0):}~>{}}},"B2"*{i},"B9"*{j},
"A1"*{\bullet},"A2"*{\bullet},"A3"*{\bullet},
"A4"*{\bullet},"A5"*{\bullet},"A6"*{\bullet},"A7"*{\bullet},
"A8"*{\bullet},"A9"*{\bullet},"A10"*{\bullet},"A2";"A9"**@{-}
\endxy
\quad \longmapsto \quad 
\xy/r2.8pc/:,{\xypolygon10"A"{~>{}{\bullet}}},
{\xypolygon10"B"{~:{(1.25,0):}~>{}}},"B2"*{i},"B9"*{j},
\POS"A2" \ar @/_1ex/ "A9" \POS"A9" \ar @/_1ex/ "A2"
\endxy$$
(All pictures here are with $n=10$.)
One easily checks that the reflecting hyperplane $H_{i,j}$ of $(i\; j)$ is
visible from $v$, and that $\bs_{i,j}$ is the corresponding local generator.
More generally, to any non-crossing partition of $\mu_n$
(cf. \cite{BDM} or \cite{reiner}),
we associate an element of $\pi_1(X_n,\mu_n)$ in the following manner:
$$
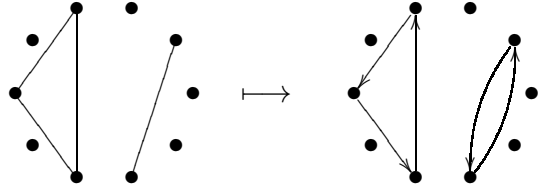
\xy/r2.8pc/:,{\xypolygon10"A"{~>{}}},
"A1"*{\bullet},"A2"*{\bullet},"A3"*{\bullet},
"A4"*{\bullet},"A5"*{\bullet},"A6"*{\bullet},"A7"*{\bullet},
"A8"*{\bullet},"A9"*{\bullet},"A10"*{\bullet},"A2";"A9"**@{-},
"A4";"A6"**@{-},"A6";"A8"**@{-},"A8";"A4"**@{-},
\endxy
\quad \longmapsto \quad 
\xy/r2.8pc/:,{\xypolygon10"A"{~>{}{\bullet}}},
{\xypolygon10"B"{~:{(1.25,0):}~>{}}},
\POS"A2" \ar @/_1ex/ "A9" \POS"A9" \ar @/_1ex/ "A2"
\POS"A4" \ar "A6" \POS"A6" \ar "A8"
\POS"A8" \ar "A4"
\endxy$$
The planar oriented graph above may be interpreted, upon need, in
three distinct but consistent ways:
as an actual path (see section 4 in \cite{BDM}),
as an element of $\pi_1(X_n,\mu_n)$,
or as the graph of the corresponding permutation,
via the morphism $\pi_1(X_n,\mu_n) \rightarrow \mathfrak{S}_n$.
The elements of $\mathfrak{S}_n$ obtained this way are exactly the 
elements of $P_c$. This correspondence is a poset isomorphism between
the poset of non-crossing partitions (for the ``is finer than'' order) and
$(P_c,\prec_T)$. Since non-crossing partitions form a lattice,
this proves Fact \ref{factB} in this case.

The Coxeter element $c$ correspond to the partition with only one
part:
$$
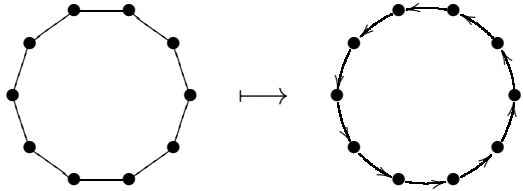
\xy/r2.8pc/:,{\xypolygon10"A"{~>{-}}},
"A1"*{\bullet},"A2"*{\bullet},"A3"*{\bullet},
"A4"*{\bullet},"A5"*{\bullet},"A6"*{\bullet},"A7"*{\bullet},
"A8"*{\bullet},"A9"*{\bullet},"A10"*{\bullet}
\endxy
\quad \longmapsto \quad
\xy/r2.8pc/:,{\xypolygon10"A"{~>{}{\bullet}}},
\POS"A1" \ar @/_0.251ex/ "A2" \POS"A2" \ar @/_0.251ex/ "A3",
\POS"A3" \ar @/_0.251ex/ "A4" \POS"A4" \ar @/_0.251ex/ "A5",
\POS"A5" \ar @/_0.251ex/ "A6" \POS"A6" \ar @/_0.251ex/ "A7",
\POS"A7" \ar @/_0.251ex/ "A8" \POS"A8" \ar @/_0.251ex/ "A9",
\POS"A9" \ar @/_0.251ex/ "A10" \POS"A10" \ar @/_0.251ex/ "A1",
\endxy$$
The Coxeter element $c$ is 
the element $c_{L,R}$, where $(L,R)$ is the
following chromatic pair (or any other chromatic pair obtained by rotating
the picture):
$$
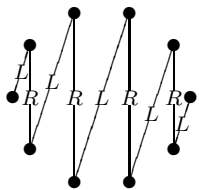
\xy/r2.8pc/:,{\xypolygon10"A"{~>{}}},
"A1"*{\bullet},"A2"*{\bullet},"A3"*{\bullet},
"A4"*{\bullet},"A5"*{\bullet},"A6"*{\bullet},"A7"*{\bullet},
"A8"*{\bullet},"A9"*{\bullet},"A10"*{\bullet},
\POS"A6" \ar @{-} "A5" |L
\POS"A5" \ar @{-} "A7" |R
\POS"A7" \ar @{-} "A4" |L
\POS"A4" \ar @{-} "A8" |R
\POS"A8" \ar @{-} "A3" |L
\POS"A3" \ar @{-} "A9" |R
\POS"A9" \ar @{-} "A2" |L
\POS"A2" \ar @{-} "A10" |R
\POS"A10" \ar @{-} "A1" |L
\endxy
$$

These pictures provide good illustrations of many of our results. For example,
conjugating by $c$ is the same as ``rotating pictures by one $n$-th of
a turn''. The isomorphism $\BT \simeq T$ from Lemma \ref{braidorbit} is
explained by the fact that the above graph picturing
the chromatic pair $(L,R)$ generates,
by rotation, the complete graph on $\mu_n$.
Lemma \ref{2hgon} is also easy to figure out: for any
$\zeta,\zeta'\in\mu_n$,
we have $\frac{\zeta-\zeta'}{|\zeta-\zeta'|}\in \mu_{2n}$.
The type $A$ case of Theorem \ref{conjecture} is also clear.

For $s,t\in T$, we have $s\parallel_c t$ if and only if the edges
corresponding to $s$ and $t$ have a common endpoint or no common point.
All relations claimed in Fact \ref{factA} are particular Sergiescu
relations (\cite{sergiescu}).

\subsection{Type $B$}
Let $(W,T)$ be the reflection group of type $B_n$, in its usual
monomial realization. The Coxeter number is $2n$.
It is well known that the orbit space
$W \backslash V^{\reg}$ can be identified with the space of subsets
of $\BC^*$ of cardinal $n$ or, equivalently, with the fixed subspace
$X_{2n}^{\mu_2}$ for the action of $\mu_2$ on $X_{2n}$.
A particular case of Proposition 5.1 in \cite{BDM} identifies
$\pi_1(X_{2n}^{\mu_2},\mu_{2n})$ with $\pi_1(X_{2n},\mu_{2n})^{\mu_2}$.

Let $x\in X_{2n}^{\mu_2}$. The identification $W \backslash V^{\reg}$ is such
that the fiber in $V^{\reg}$ above $x$ is the set of $n$-tuples
$(x_1,\dots,x_n)\in\BC^n$ such that $x=\{x_1,\dots,x_n,-x_1,\dots,-x_n\}$.
In particular, the vector $$v:=(e^{2i\pi\frac{1}{2n}},e^{2i\pi\frac{2}{2n}},
\dots,e^{2i\pi\frac{n}{2n}})$$ lies above $\mu_{2n}$.
It is a regular $e^{2i\frac{\pi}{2n}}$-eigenvector for the Coxeter element
 $$c:= \mbox{ \tiny $\begin{pmatrix}
0 & 1 & 0 & \hdotsfor{2} & 0 \\
0 & 0 & 1 & 0 & \dots & 0 \\
\hdotsfor{6} \\
\hdotsfor{6} \\
0 & \hdotsfor{3} & 0 & 1 \\
-1 & 0 & \hdotsfor{3} & 0 
\end{pmatrix}$ }$$

We say that a partition of $\mu_{2n}$ is $\mu_2$-symmetric if
each part is stable by multplication by $-1$.
Let $\lambda$ be a $\mu_2$-symmetric
non-crossing partition of $\mu_{2n}$.
Let $\sigma_{\lambda}$ be the corresponding element
of $\mathfrak{S}_{2n}$ (identified, as in the type $A$ discussion
above, with $\mathfrak{S}_{\mu_{2n}}$).
For any $k\in\{1,\dots,n\}$, there is unique pair
$(l_k,\varepsilon_k)\in \{1,\dots,n\} \times \{\pm 1\}$ such that
$\sigma_{\lambda}(e^{2i\pi\frac{k}{2n}})=
\varepsilon_k e^{2i\pi\frac{l_k}{2n}}$. To $\lambda$, we associate
the monomial matrix $w_{\lambda} := (\varepsilon_k \delta_{l_k,l})_{l,k}$ in
$W$. An example with $n=5$ is illustrated below:

$$\lambda: 
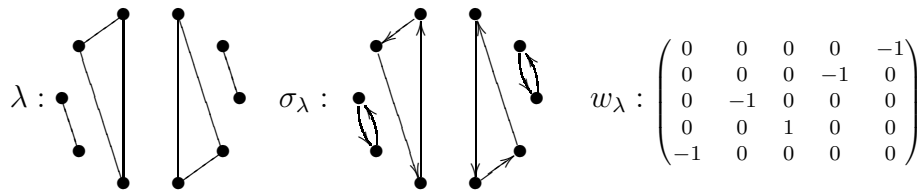
\xy/r2.8pc/:,{\xypolygon10"A"{~>{}}},
"A1"*{\bullet},"A2"*{\bullet},"A3"*{\bullet},
"A4"*{\bullet},"A5"*{\bullet},"A6"*{\bullet},"A7"*{\bullet},
"A8"*{\bullet},"A9"*{\bullet},"A10"*{\bullet},
"A4";"A5"**@{-},"A5";"A8"**@{-},"A8";"A4"**@{-},
"A6";"A7"**@{-},"A1";"A2"**@{-},
"A9";"A10"**@{-},"A10";"A3"**@{-},"A3";"A9"**@{-},
\endxy
\quad
\sigma_{\lambda}: \xy/r2.8pc/:,{\xypolygon10"A"{~>{}{\bullet}}},
{\xypolygon10"B"{~:{(1.25,0):}~>{}}},
\POS"A6" \ar @/_0.5ex/ "A7" \POS"A7" \ar @/_0.5ex/ "A6"
\POS"A1" \ar @/_0.5ex/ "A2" \POS"A2" \ar @/_0.5ex/ "A1"
\POS"A4" \ar "A5" \POS"A5" \ar "A8" \POS"A8" \ar "A4"
\POS"A9" \ar "A10" \POS"A10" \ar "A3" \POS"A3" \ar "A9"
\endxy
\quad
w_{\lambda}:
\mbox{\tiny $\begin{pmatrix}
0 & 0 & 0 & 0 & -1 \\
0 & 0 & 0 & -1 & 0 \\
0 & -1 & 0 & 0 & 0 \\
0 & 0 & 1 & 0 & 0 \\
-1 & 0 & 0 & 0 & 0
\end{pmatrix}$ }
$$
One can easily deduce from the type $A$ case that this construction
identifies the poset of $\mu_2$-symmetric non-crossing partitions
of $\mu_{2n}$ with $P_c$. The type $B$ case of 
Fact \ref{factB} follows, since
$\mu_2$-symmetric partitions form a lattice (this lattice is studied
in \cite{reiner}).

The reflections in $W$ correspond to minimal symmetric non-crossing
partitions. There are two types of them, corresponding to the two
conjugacy classes of reflections in $W$: partitions with
one symmetric part
$\{\zeta,-\zeta\}$ (and all
other parts being points), and partitions with
two non-symmetric parts $\{\zeta,\zeta'\}$
and $\{-\zeta,-\zeta'\}$ (with $\zeta\neq \pm \zeta'$), as illustrated
below. We call the first type ``long'' and the second ``short''.

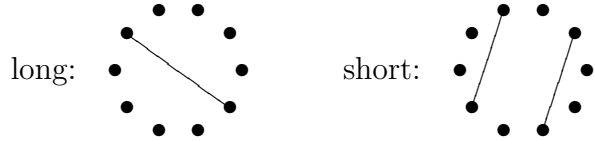
$$\text{long:} \quad \xy/r2pc/:,{\xypolygon10"A"{~>{}}},
"A5";"A10"**@{-},
"A1"*{\bullet},"A2"*{\bullet},"A3"*{\bullet},
"A4"*{\bullet},"A5"*{\bullet},"A6"*{\bullet},"A7"*{\bullet},
"A8"*{\bullet},"A9"*{\bullet},"A10"*{\bullet}
\endxy
\quad \quad \quad
\text{short:} \quad \xy/r2pc/:,{\xypolygon10"A"{~>{}}},
"A2";"A9"**@{-},"A4";"A7"**@{-},
"A1"*{\bullet},"A2"*{\bullet},"A3"*{\bullet},
"A4"*{\bullet},"A5"*{\bullet},"A6"*{\bullet},"A7"*{\bullet},
"A8"*{\bullet},"A9"*{\bullet},"A10"*{\bullet}
\endxy$$
Here again, the corresponding braids are the local generators at $v$,
two reflections are non-crossing if and only if the corresponding
edges have no common point (except possibly endpoints; the two
reflections pictured above are crossing)
and Fact \ref{mainA} follows from the usual type $A$ Sergiescu relations.

\subsection{Type $D$}
Let $(W,T)$ be the reflection group of type $D_n$, with $n\geq 3$, seen
in its usual monomial realization. The degrees of $D_n$ are 
$2,4,6,\dots,2(n-1),n$. The Coxeter number is $2(n-1)$.
A Coxeter element is
 $$c:= \mbox{\tiny $\begin{pmatrix}
0 & 1 & 0 & \hdotsfor{2} & 0 & 0\\
0 & 0 & 1 & 0 & \dots & 0  & : \\
\hdotsfor{6} & :\\
\hdotsfor{6} & :\\
0 & \hdotsfor{3} & 0 & 1 & :\\
-1 & 0 & \hdotsfor{3} & 0 & 0 \\
0 & \hdotsfor{4} & 0 & -1
\end{pmatrix} $ }$$
(the matrix has two diagonal blocks: a $(n-1)\times (n-1)$ block 
corresponding to a type $B_{n-1}$ Coxeter element, and $-1$ as 
last diagonal coefficient). As this matrix suggests, the
dual geometry of type $D_n$ is related to the dual geometry
of type $B_{n-1}$.

A regular 
$e^{2i\frac{\pi}{2(n-1)}}$-eigenvector for $c$ is
$$v:=(e^{2i\pi\frac{1}{2(n-1)}},e^{2i\pi\frac{2}{2(n-1)}},
\dots,e^{2i\pi\frac{n-1}{2(n-1)}},0).$$

Consider the map 
\begin{eqnarray*}
p:V^{\reg} & \rightarrow & \mathcal{P}(\BC) \\
(x_1,\dots,x_n) & \mapsto & \{x_1,\dots,x_n,-x_1,\dots,-x_n\}
\end{eqnarray*}
For $1\leq i,j \leq n$, we denote $H_{i,j}$ (resp. $\overline{H}_{i,j}$)
the reflecting
hyperplane with equation $X_i=X_j$ (resp. $X_i=-X_j$).
We denote by $s_{i,j}$ and $\overline{s}_{i,j}$
the corresponding reflections. Contrary to the
type $B$ case, the hyperplanes with equation $X_i=0$ are not reflecting
hyperplanes. The image of $p$ is in $X_{2n}\cup X_{2n-1}$. We have
$p(v) = \mu_{2(n-1)} \cup \{0\} \in X_{2n-1}$.

For $i\in\{1,\dots,n-1\}$, we set $\zeta_i:=e^{2i\pi \frac{i}{2(n-1)}}$.
 If $1\leq i<j \leq n-1$, we represent the
reflection $s_{i,j}$ (resp. $\overline{s}_{i,j}$)
by the planar graph on $p(v)$ with edges
$[\zeta_i,\zeta_j]$ and $[\overline{\zeta}_i,\overline{\zeta}_j]$ 
(resp. $[\zeta_i,\overline{\zeta}_j]$ and $[\overline{\zeta}_i,{\zeta}_j]$).
If $1\leq i \leq n-1$, we represent the reflection $s_{i,n}$ 
(resp. $\overline{s}_{i,n}$) by the
planar graph with only edge $[\zeta_i,0]$ (resp. $[\overline{\zeta}_i,0]$).
Here are some examples with $n=6$:

$$s_{3,4}: 
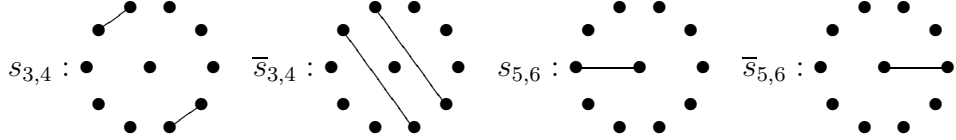
\xy/r2pc/:,{\xypolygon10"A"{~>{}}},
"A4";"A5"**@{-},"A9";"A10"**@{-},
"A0"*{\bullet},"A1"*{\bullet},"A2"*{\bullet},"A3"*{\bullet},
"A4"*{\bullet},"A5"*{\bullet},"A6"*{\bullet},"A7"*{\bullet},
"A8"*{\bullet},"A9"*{\bullet},"A10"*{\bullet}
\endxy \quad
\overline{s}_{3,4}: \xy/r2pc/:,{\xypolygon10"A"{~>{}}},
"A4";"A10"**@{-},"A9";"A5"**@{-},
"A0"*{\bullet},"A1"*{\bullet},"A2"*{\bullet},"A3"*{\bullet},
"A4"*{\bullet},"A5"*{\bullet},"A6"*{\bullet},"A7"*{\bullet},
"A8"*{\bullet},"A9"*{\bullet},"A10"*{\bullet}
\endxy \quad
s_{5,6}: \xy/r2pc/:,{\xypolygon10"A"{~>{}}},
"A0";"A6"**@{-},
"A0"*{\bullet},"A1"*{\bullet},"A2"*{\bullet},"A3"*{\bullet},
"A4"*{\bullet},"A5"*{\bullet},"A6"*{\bullet},"A7"*{\bullet},
"A8"*{\bullet},"A9"*{\bullet},"A10"*{\bullet}
\endxy \quad
\overline{s}_{5,6}: \xy/r2pc/:,{\xypolygon10"A"{~>{}}},
"A0";"A1"**@{-},
"A0"*{\bullet},"A1"*{\bullet},"A2"*{\bullet},"A3"*{\bullet},
"A4"*{\bullet},"A5"*{\bullet},"A6"*{\bullet},"A7"*{\bullet},
"A8"*{\bullet},"A9"*{\bullet},"A10"*{\bullet}
\endxy$$
We say that $s_{i,j}$ (resp. $\overline{s}_{i,j}$) is $B$-like if
both $i<n$ and $j<n$. This notion is of course specific to our choice
of $c$.

We leave to the reader the following lemma:
\begin{lemma}
Two reflections in $T$ are non-crossing (with respect to $c$) if
and only if the associated graphs are non-crossing (i.e, their
edges have no common points except possibly endpoints).
\end{lemma}

The corresponding local generators at $v$ are easy to compute.
We may represent each of them by a path in $V^{\reg}$
starting at $v$ and ending at $sv$, according to the pictures below:
$$\bs_{3,4}: 
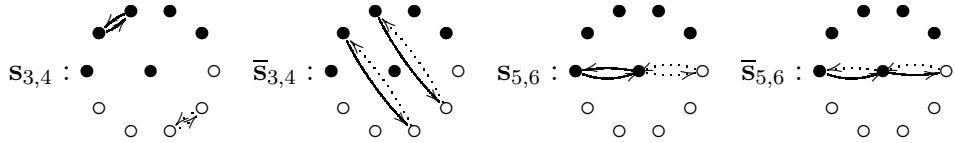
\xy/r2pc/:,{\xypolygon10"A"{~>{}{\circ}}},
"A0"*{\bullet},"A2"*{\bullet},"A3"*{\bullet},
"A4"*{\bullet},"A5"*{\bullet},"A6"*{\bullet},
\POS"A4" \ar @/_0.25ex/ "A5" \POS"A5" \ar @/_0.25ex/ "A4"
\POS"A9" \ar @{.>}@/_0.25ex/ "A10" \POS"A10" \ar @{.>}@/_0.25ex/ "A9"
\endxy \quad
\overline{\bs}_{3,4}: \xy/r2pc/:,{\xypolygon10"A"{~>{}{\circ}}},
"A0"*{\bullet},"A2"*{\bullet},"A3"*{\bullet},
"A4"*{\bullet},"A5"*{\bullet},"A6"*{\bullet},
\POS"A4" \ar @/_0.4ex/ "A10" \POS"A5" \ar @/_0.4ex/ "A9"
\POS"A9" \ar @{.>}@/_0.4ex/ "A5" \POS"A10" \ar @{.>}@/_0.4ex/ "A4"
\endxy \quad
\bs_{5,6}: \xy/r2pc/:,{\xypolygon10"A"{~>{}{\circ}}},
"A0"*{\bullet},"A2"*{\bullet},"A3"*{\bullet},
"A4"*{\bullet},"A5"*{\bullet},"A6"*{\bullet},
\POS"A6" \ar @/_0.5ex/ "A0" \POS"A0" \ar @/_0.25ex/ "A6"
\POS"A1" \ar @{.>}@/_0.5ex/ "A0" \POS"A0" \ar @{.>}@/_0.25ex/ "A1"
\endxy \quad
\overline{\bs}_{5,6}: \xy/r2pc/:,{\xypolygon10"A"{~>{}{\circ}}},
"A0"*{\bullet},"A2"*{\bullet},"A3"*{\bullet},
"A4"*{\bullet},"A5"*{\bullet},"A6"*{\bullet},
\POS"A6" \ar @/_0.5ex/ "A0" \POS"A0" \ar @/_0.25ex/ "A1"
\POS"A1" \ar @{.>}@/_0.5ex/ "A0" \POS"A0" \ar @{.>}@/_0.25ex/ "A6"
\endxy$$
These pictures should be interpreted as follows: the black dots indicates
the starting values of the coordinates (the coordinates of $v$); these
coordinates vary continuously according to the plain arrows. The white
dots and the dotted arrows complete the picture by symmetry. Together,
the plain and dotted arrows represent the image of the path by $p$.

\smallskip
{\bf \flushleft Remark.}
The ``folding'' of the $D_n$ Dynkin diagram onto the $B_{n-1}$ diagram
has a dual analog: the type $B_{n-1}$ dual monoid is isomorphic to the
submonoid of the type $D_n$ dual monoid generated by the (short)
$B$-like
reflections and the (long)
products $s_{k,n}\overline{s}_{k,n} (= \overline{s}_{k,n}
s_{k,n})$, for $k=1,\dots,n-1$.

\smallskip

The dual relations needed for Fact \ref{mainA} are easy variations on 
Sergiescu relations, left to the reader.
Using Allcock's ``orbifold'' pictures (which are quotient modulo $\pm 1$ of the pictures used here), Picantin explicitly described a presentation of
the type $D$ braid monoid (\cite{picantin}); Fact \ref{mainA} can also
be checked in the presentation in \cite{picantin}.

A detailed combinatorial proof of Fact \ref{mainB} in type $D$ can be found
in \cite{bradywatt}. Let us sketch a more geometric proof.
As noted after Fact \ref{mainB}, the lattice property
would follow if we prove that any pair of crossing reflections has
a right lcm. The case of two $B$-like reflections follows easily from the
type $B$ combinatorics. Since any pair of crossing reflections always
contain a $B$-like reflection, the only case left is when a $B$-like
reflection $s_{i,j}$ (or $\overline{s}_{i,j}$)
is crossing with a reflection of the form $s_{k,n}$.
The lcm may be computed explicitly, using convex hulls.

Brady-Watt and Picantin note that, in types $A$ and $B$, the lattice $P_c$
is isomorphic to the corresponding 
lattice of non-crossing partitions, defined by Reiner. The local geometry
at a $h$-regular eigenvector provides an explanation.
In type $D$, it ought to be possible to encode elements of $P_c$ 
by planar graphs (extending what is done here for reflections), such
that this encoding provides an isomorphism
between Reiner's type $D$ lattice of
non-crossing partitions and $P_c$ (see also \ref{catalan} below).

\section{Numerology}

\subsection{The duality}
\label{duality}
This subsection is an attempt to convince the reader that there is some
sort of ``duality'' between the classical braid monoid and the dual
braid monoid.
Unfortunately, at the present time, we are not
able to formalize the nature of this duality.

Let $(W,S=L\cup R)$ be an irreducible finite Coxeter system of rank $n$, with
set of reflections $T$. The notation $N:=|T|$ is standard.
The set of atoms of $\BB_+(W,S)$ is $\BS$. 
Let $\bc:=\bc_{L,R}$. We denote by $p$ the morphism
$\BB(W,S)\rightarrow W, \bs \mapsto s$.
Let $c:=p(\bc)=c_{L,R}$.
We use Corollary
\ref{coroB} to identify $\BM(P_c)$ with the submonoid of $\BB(W,S)$
generated by $\BT = \cup_{k\in \BZ} \bc^k \BS \bc^{-k}$.

The lcm (in $\BB_+(W,S)$)
of the atoms of $\BB_+(W,S)$ is $\bw_0$; it has length
$N$ for the natural length function on $\BB(W,S)$; it image
$w_0$ in $W$ has order $2$.
The lcm (in $\BM(P_c)$)
of the atoms of $\BM(P_c)$ is $\bc$; it has length $N$; the 
order of $c$ is the Coxeter number $h$.

Write $L=\{s_1,\dots,s_k\}$, $R=\{s_{k+1},\dots,s_n\}$.
We have $\bc=\bs_1\dots \bs_n$. In other words, $\bc$ is the product
of the atoms of $\BB_+(W,S)$, taken in a suitable order.
Similarly, the next lemma proves that
the product of atoms of $\BM(P_c)$, taken in a suitable order,
is $\bw_0$.

We extend the notation $\bs_m$ to all positive integers $m$, in such
a way that $\bs_m$ only depends on $m\mod n$.

\begin{lemma}
For any positive integer $m$, set 
$$\bt_m := \left( \prod_{i=1}^m \bs_i \right) \left( \prod_{i=1}^{m-1} \bs_i
\right)^{-1}.$$
We have $\BT = \{\bt_1,\dots,\bt_N\}$, and
$$\bw_0 = \prod_{m=1}^{N} \bt_{N-m+1}.$$
\end{lemma}

\begin{proof}
Set $t_m:=p(\bt_m)$.
By \cite{bourbaki}, Ch. V, \S 6 Ex. 2, pp 139-140,
we have $T = \{t_1,\dots,t_N\}$.

Using the commutation relations within $L$, we see that, when
$1\leq m \leq k$, $\bt_m = \bs_m$.
When $k+1 \leq m \leq n$, the commutation relations within
$R$ yield
$$\bt_m = \bs_1 \dots \bs_{m-1} \bs_m 
\bs_{m-1}^{-1} \dots \bs_1
= \bs_1 \dots \bs_n \bs_m \bs_n^{-1} \dots \bs_1
= \bc \bs_m \bc^{-1}.$$
We have proved
$$\{\bs_1,\dots,\bs_k,\bc \bs_{k+1} \bc^{-1},\dots, \bc \bs_n \bc^{-1} \}
= \{\bt_1,\dots,\bt_n\}.$$

For all $m$, we have $\bt_{m+n} = \bc \bt_m \bc^{-1}$.
From this and the above description of $\{\bt_1,\dots,\bt_n\}$, we deduce
$\BT = \{\bt_1,\dots,\bt_N\}$.

From the Bourbaki exercice quoted above, we also get
$w_0=\prod_{m=1}^N s_m$. Since
$(s_1,\dots,s_N) \in \Red_S(w_0)$, we have
$\bw_0 = \prod_{m=1}^N \bs_m = \prod_{m=1}^N \bt_{N-m+1}$.
\end{proof}

These facts are summarized in the following table:

%\begin{table}
\begin{center}
\begin{tabular}{|c|c|c|}
\hline
& Classical monoid & Dual monoid \\
%\hline
%Length function & $l_S$ & $l_T$ \\
\hline
Set of atoms &  $\BS$ & $\BT$ \\
Number of atoms & $n$ & $N$ \\
\hline
$\Delta$ & $\bw_0$ & $ {\bc}$ \\
Length of $\Delta$ & $N$ & $n$ \\
Order of $p(\Delta)$ & $2$ & $h$ \\
\hline
Product of the atoms & $\bc$ & $\bw_0$ \\
Regular degree & $h$ & $2$ \\
\hline
\end{tabular}
\end{center}
%\end{table}

The final line has the following explanation: 
in \cite{zariski}, a certain class of presentations of 
braid groups is constructed. Each of these presentations 
corresponds to a regular degree $d$.
The product of the generators, raised to the
power $d$ (which is the order of the image of this product
in the reflection group), is always 
central. 

For an irreducible Coxeter group, $2$ and $h$ are the respectively
smallest and
largest
degrees; they are always regular;
it is possible to choose intermediate regular degrees but they do
not seem to yield Garside monoids.

\subsection{Catalan numbers}
\label{catalan}
Reiner (\cite{reiner}) suggests a definition
for what should be the ``Catalan number'' attached to a finite Coxeter group.
Though he has no definition for exceptional types non-crossing partitions,
the Catalan number should be the number of non-crossing partitions.
The usual Catalan numbers
$$c_n:=\frac{1}{n+1}\begin{pmatrix} 2n \\ n \end{pmatrix}$$
correspond to type $A$.

For type $A$ and $B$,
Brady-Watt and Picantin (\cite{bradywatt}, \cite{picantin})
noticed the relation between Reiner's non-crossing partitions and
the dual monoid.

At the level of Catalan numbers, the correspondance works for all types:

\begin{prop}
Let $W$ be an irreducible Coxeter group, with degrees $d_1,\dots,d_n=h$.
Let $c$ be a Coxeter element. The
number of simple elements in the dual monoid is given by the formula
$$| P_c |= \prod_{i=1}^n \frac{d_i + h}{d_i}.$$
\end{prop}

\begin{proof}
This can be checked, case-by-case, on the list given in
\cite{picantin}
\end{proof}

We may now answer some of
the questions raised in \cite{reiner}, Remark 2.
Fix a Coxeter element $c$. An element $w\in W$ should be called
{\sf non-crossing} if it is in $P_c$. A subspace in the intersection
lattice generated by the reflecting hyperplanes should be called
{\sf non-crossing} if it is of the form $K_w$, with $w$ non-crossing.
According to Proposition
\ref{brady}, this defines a one-to-one correspondence
between non-crossing
elements and non-crossing subspaces (and standard parabolic
subgroups). A more detailed study of the local
geometry at a $h$-regular eigenvector is likely to provide explanation

\begin{question}
The function $l_T$ gives a natural grading on $P_c$. What should
be the formula for the Poincar\'e
polynomial of $P_c$?
\end{question}

{\small
{\flushleft \bf Example.}
For the reflection group of type $E_8$, this Poincar\'e polynomial
is 
$$1 + 120q + 1540 q^2+6120 q^3 + 9518 q^4 + 6120 q^5+1540 q^6 + 120 q^7 + q^8$$
(the palindromicity of this polynomial is a general fact,
easy to prove: consider the
bijection $P_c \rightarrow P_c$, $w\mapsto w^{-1} c$).
The value of this polynomial at $q=1$ is the cardinal of $P_c$ (the
corresponding Catalan number). Here, this value is $25080$.
Note that the order of $W(E_8)$ is
$696729600$; while the enumeration of the elements of $W(E_8)$ is
presently beyond reach, the poset $P_c$ is small enough to be enumerated 
by computer; 
checking the lattice property does not require much computing power.
Using Lemma \ref{carter} and the classical formula of Solomon,
we see that the Poincar\'e polynomial of $W(E_8)$ for the length function
$l_T$ is
{ $$(1+q)(1+7q)(1+11q)(1+13q)(1+17q)(1+19q)(1+23q)(1+29q)=$$
$$1 + 120 q + 6020 q^2 + 163800 q^3 + 2616558q^4
+ 24693480 q^5 + 130085780 q^6 $$ $$ + 323507400 q^7 + 215656441 q^8.$$}
A final remark about the arithmetic of $P_c$. In type $A_n$, the lcm 
of two crossing reflections has length $3$. In type $E_8$, there
are pairs of crossing reflections whose lcm is $c$, of length $8$.
}

\section{Applications and problems}
\label{applications}

\subsection{The dual normal form}

As mentioned in the appendix, Garside monoids admit natural normal
forms. Therefore, the dual monoid yields a new solution to the word 
problem. In the type $A$, the complexity of this solution has been 
studied by Birman-Ko-Lee (\cite{BKL}) and has been proved to be better
than the one deriving from the classical monoid.
The general case has yet to be studied. A possible advantage that can
already be observed is that the Catalan number $|P_c|$ is much smaller
than $|W|$ (in the $E_8$ example above, $|P_c|$ is not far from being
the square root of $|W|$).

Another specificity of the dual normal form 
is that it is compatible with the conjugacy action
of a Coxeter element. In \cite{BDM} is mentioned a conjecture about
centralizers in generalized 
braid groups of certain $d$-th roots of central elements,
and the Birman-Ko-Lee monoid is used to prove the conjecture for
the type $A$ case.
A specificity of the Birman-Ko-Lee monoid, used in section 4 of \cite{BDM},
is that it is possible to associate to each element $P_c$ a ``geometric
normal form'' (a particular loop which is the  shortest 
loop in its homotopy class, for a suitable metric).
We suspect the same can be done with the dual braid monoid.
New cases for the centralizer conjecture would follow (the case of $W$ being 
a Coxeter group, and $d$ dividing the Coxeter number $h$).

\subsection{Braid groups actions on categories}
\label{actions}

To illustrate how the dual monoid can be used as a replacement for the classical
braid monoid, we discuss the problem of braid groups actions on categories.
This problem has been studied by Deligne; the present
discussion is nothing more than a straightforward reformulation of
\cite{deligne2} in the more general context of Garside monoids.

An action of a pre-monoid $P$ on a category $\CC$ is a collection
of endofunctors $(T(f))_{f\in P}$ and of natural isomorphisms
$c_{f,g}: T(f) \circ T(g) \rightarrow T(fg)$ (one for each pair $(f,g)$ in
the domain of the partial product) with the following compatibility condition:
whenever the product $fgh$ is defined in $P$, the diagram 
$$\xymatrix{T(f) \circ T(g) \circ T(h) \ar[r] \ar[d] & T(fg) \circ T(h) \ar[d]\\
T(f) \circ T(gh) \ar[r] & T(fgh) }$$
is commutative.

For monoids, one recovers the notion of action 
on a category defined in \cite{deligne2}.
An action of $\BM(P)$ on $\CC$ gives, by restriction, an action of $P$ on
$\CC$. The actions of $\BM(P)$ (resp. $P$) on $\CC$ form a category and
the restriction is functorial.
The analog for the dual braid monoid of the main result (Th\'eor\`eme 1.5)
in \cite{deligne2} is a special case of:

\begin{theo}[after Deligne]
Let $P$ be a Garside pre-monoid. Let $\CC$ be a category.
The restriction functor from the category of actions of $\BM(P)$ on $\CC$ to
the category of actions of $P$ on $\CC$ is an equivalence of categories.
\end{theo}

(An unital action by auto-equivalences of $\BM(P)$ extends to an action 
of the braid group $\BG(P)$; see \cite{deligne2}, Prop. 1.9.)

The construction of the quasi-inverse is virtually identical to the one
in \cite{deligne2}, and the proofs can be reproduced with only minor
adaptations. 
Let $m\in \BM(P)$. We define in the appendix a poset $(E(m),\leq)$
(our definition actually mimics the one from {\em loc. cit.}).
Generalizing Th\'eor\`eme 2.4 in {\em loc. cit.}, one can prove that
the geometric realization $|E(m)|$ is contractible. 

Denote by $U$ the set of atoms of $P$ which are left-divisors of $m$.
For $u\in U$, denote by $E_u(m)$ the subset of $E(m)$
consisting of those sequences 
$(p_1,\dots,p_k)$ such that $u\prec p_1$. For any 
non-empty subset $V\subset U$, let
$E_V(m):=\cap_{u\in V} E_u(m)$.
Deligne's proof can be easily adapted to establish
the contractibility of $|E(m)|$.
For the convenience of the reader, we precise how
Lemme 2.5 and its proof should be modified to get rid of galleries and
chambers:

\begin{lemma}
Let $m\in \BM(P)$. Let $U$ be as above. Let $V$ be a non-empty subset of $U$.
The geometric realization $|E_V(m)|$ is 
contractible.
\end{lemma}

\begin{proof}
Let $\delta_V$ be the right lcm of the elements of $V$.
Since $\forall u\in V, u\prec m$, one has $\delta_V\prec m$. More
precisely, for all $(p_1,\dots,p_k)\in E_V(m)$, one has $\delta_V\prec p_1$.
Let $n\in M(P)$ be the element uniquely defined by $\delta_V n = m$.
As $V$ is non-empty, one has $l(n)<l(m)$, and Deligne's proof's induction
hypothesis implies that $|E(n)|$ is contractible.

The map
\begin{eqnarray*}
f: E(n) & \longrightarrow & E_V(m) \\
(p_1,\dots,p_k) & \longmapsto & (\delta_V,p_1,\dots,p_k)
\end{eqnarray*}
is increasing and induces an isomorphism between 
$E(n)$ and an initial segment of $E_V(m)$.

The map
\begin{eqnarray*}
f^*: E_V(m) & \longrightarrow & E(n) \\
(p_1,\dots,p_k) & \longmapsto & 
\begin{cases}
(\delta_V^{-1}p_1,\dots,p_k) & \text{if $\delta_V\neq p_1$}\\
(p_2,\dots,p_k) & \text{if $\delta_V = p_1$}
\end{cases}
\end{eqnarray*}
is increasing and one has, for all $x\in E(n)$ and all $y\in E_V(m)$,
$$f(x) \leq y \Leftrightarrow x \leq f^*(y),$$
and one concludes as in Deligne's proof.
\end{proof}

\subsection{New $K(\pi,1)$'s for braid groups}
A motivation for Brady's work on the Birman-Ko-Lee monoid was
to construct new
finite simplicial complexes which are $K(\pi,1)$'s for braid groups
(\cite{brady}, \cite{bradywatt}).
His techniques are modelled on a construction of Bestvina.
Following the same approach,
Charney, Meier and Whittlesey have extended Bestvina's construction
to the context of an arbitrary Garside monoid (\cite{charney}).
They note that the $K(\pi,1)$ constructed from the dual monoid
has the minimal possible dimension.

For a general Garside monoid, the $K(\pi,1)$ constructed in \cite{charney}
is related to the complexes $E(m)$ from the
previous subsection (more specifically, to $E(\Delta)$, where $\Delta$
is the Garside element).

\subsection{Problems}
We conclude with a list of problems.

(1) Formalize and complete the ``dual Coxeter theory'' of the first section.
What should be the ``dual exchange lemma''?

(2) What can be done with infinite Coxeter groups?

(3)
Provide proofs of Facts \ref{factA} and \ref{factB} not relying on the 
classification of finite Coxeter systems.

(4) Classify all local monoids which are Garside.

(5) Study the relations between the 
three natural orders on $W$: $\prec_S$, $\prec_T$, and the Bruhat order.
Does the order $\prec_T$ have a 
geometric interpretation similar to the ones known for the Bruhat order?

(6) Study Hecke algebras with the dual point of view. Elements
$T_w$ are easy to define when $w\in P_c$. The work of
Bremke-Malle is a possible source of inspiration on how to define
$T_w$ when $w\notin P_c$ (see \cite{bremkemalle},
Prop. 2.4). 
More generally,
study objects classified by Weyl groups (Lie groups, algebraic
groups,...) with the dual point of view.

(7) Explain and formalize
the ``duality'' between the classical and the dual monoid.

(8) There should be a bijection between $P_c$
and the number of regions inside the fundamental chamber in
the affine hyperplane arrangement described p. 219 in \cite{reiner}.
Give a general construction of such a bijection.

(9) (Related to (4) and (7)) Let $M$ be a Garside monoid.
Is it a frequent phenomenon to have another Garside monoid $N$ such
that $\BG(M)\simeq \BG(N)$?
The pair classical monoid/dual monoid is an example.
Here is another one:
the fundamental group $T_{m,n}$ of the complement of the torus link $L_{m,n}$
(obtained by closing on itself the type $A_{n-1}$ braid
$(\sigma_1\dots\sigma_{n-1})^m$)
has the presentation with $m$ generators $s_1,\dots,s_m$ and relations
$$\underbrace{s_1s_2s_3\dots}_{\text{$n$ terms}} = 
\underbrace{s_2s_3s_4\dots}_{\text{$n$ terms}}  = \dots
= \underbrace{s_ms_1s_2\dots}_{\text{$n$ terms}} $$
(if $n>m$, the $s_i$ are cyclically repeated). As noted in \cite{depa} (section
5, example 5),
the monoid $M_{m,n}$ defined by this positive presentation is a Garside
monoid.
But the links $L_{m,n}$ and $L_{n,m}$ are isotopic.
So $\BG(M_{m,n}) \simeq \BG(M_{n,m})$.
Some of these groups appear as braid groups attached to certain
complex reflection groups: according to the tables of \cite{brmarou},
$\BB(G_{12})\simeq T_{3,3}$, $\BB(G_{13})\simeq T_{3,4}$ and
$\BB(G_{22}) \simeq T_{3,5}$; hence we may define ``dual monoids''
for these braid groups.

\appendix

\section{Garsiditude}

According to the MathSciNet database, F. A. Garside published
only one mathematical paper (\cite{garside}, 1969). 
It contains a solution to the word problem in type A Artin groups.
In 1972 appeared simultaneously two articles, by Brieskorn-Saito and by
Deligne, generalizing Garside's techniques and results
to the context of an arbitrary finite type Artin group
(\cite{brieskornsaito}, \cite{deligne}).

These generalizations were however not ultimate, in the sense that most
of the crucial proofs actually work with more general groups. 
The notions of {\em Garside group} and {\em Garside monoid} 
were introduced by Dehornoy-Paris
(\cite{depa}, 1999).
A Garside group is a group which satisfies a certain number of axioms, 
sufficient to apply the techniques of Garside, Brieskorn-Saito and
Deligne. A slightly different axiomatic was introduced independently
by Corran (\cite{corran}).

Inspired by an earlier work of Michel (\cite{michel} -- which itself
reformulates results of Charney),
we proposed with Digne and Michel a variant approach to
Garside monoids (and ``locally Garside monoids''),
via what we called in \cite{BDM} {\em pre-Garside structures};
the properties of the monoid are derived from properties of a subset
of the monoid, on which the product is only partially defined.
This
approach is more natural for (classical and dual) braid monoids.
The purpose of this appendix is to give a survey of this approach, with
a new language. 
The terminology is probably more abstract than required,
not that we have any pretention to universality or exhaustivity,
but rather that we feel it simplifies the exposition.

\subsection{Pre-monoids}

A pre-monoid can be thought of as a ``fragment of monoid'' or, more 
metaphorically, as a ``seed'' containing all the information to build
a monoid.

\begin{defi}
A {\sf pre-monoid} is a triple $(P,D,m)$, where $P$ is a set,
$D$ is a subset of $P\times P$ and $m$ is a map $D\rightarrow P$, 
satisfying:
\begin{itemize}
\item[(assoc)] For all $a,b,c\in P$, the condition
``$(a,b)\in D$ and
$(m(a,b),c)\in D$'' is equivalent to ``$(b,c)\in D$ and $(a,m(b,c))\in D$'',
and, when they are satisfied, one has $m(m(a,b),c) = m(a,m(b,c))$. 
\end{itemize}
A pre-monoid $P$ is {\sf unitary} when it satisfies in addition:
\begin{itemize}
\item[(unit)] There exists an element $1\in P$, such that, for all $a\in P$,
$(a,1)\in D$ and $(1,a)\in D$, and $m(a,1)=m(1,a)=a$.
\end{itemize}
\end{defi}

The map $m$ should be seen as a ``partial product'', with domain $D$.
Practically, it is
convenient to omit to explicitly refer to $m$ and $D$: we write
``$ab$'' for ``$m(a,b)$'', and ``$ab$ is defined'' instead of
``$(a,b)\in D$''.
A trivial lemma on binary trees shows that, thanks to the
(assoc) axiom, for any
sequence $a_1$, $a_2,\dots,a_n$ of elements of $P$, the fact that the
product $a_1a_2\dots a_n$ is defined, and its value, do not depend
on how one chooses to put brackets.

Let $P$ be a pre-monoid. Let $p,q\in P$. We say that $p$ is left (resp. right)
divisor of $q$, or equivalently that $q$ is a right (resp. left)
multiple of $p$,
and we write $p \prec q$ (resp. $q\succ p$), if there exists $r\in P$ such that
$pr = q $ (resp. $rq=p$) in $P$.

\subsection{The functor $\BM$}

Pre-monoids form a category $\mathbf{preMon}$, where a morphism
$\varphi: P \rightarrow P'$ between two pre-monoids is defined to be
a set-theoretical map such that, for all $a,b\in P$ such that
$ab$ is defined, the product $\varphi(a)\varphi(b)$ is defined in $P'$,
and equal to $\varphi(ab)$.

The category $\mathbf{Mon}$ of monoids can be defined as the full subcategory
of $\mathbf{preMon}$ with objects being those pre-monoids for which
the product is everywhere defined.
The embedding functor $\mathbf{Mon} \rightarrow \mathbf{preMon}$ has
a left adjoint $\mathbf{M}$, defined as follows:
\begin{itemize}
\item Let $P$ be a pre-monoid. Let $P^*$ be the free monoid on $P$,
{\em i.e.}, the set of finite sequences of elements of $P$, for the
concatenation product. Let $\sim$ be the smallest equivalence
relation on $P^*$ compatible with concatenation and satisfying
$(a,b)\sim (ab)$ whenever $ab$ is defined in $P$. We
set $\mathbf{M}(P) := P^*/\sim$.

Note that one has a natural pre-monoid morphism
$P \rightarrow \BM(P), p \mapsto (p)$.

\item If $\varphi: P \rightarrow Q$ is a pre-monoid morphism,
we take $\BM(\varphi)$ to be the (unique) monoid morphism
which makes the following diagram commute:
$$\xymatrix{P \ar[r] \ar[d]_\varphi & \BM(P) \ar[d]^{\BM(\varphi)} \\
Q \ar[r] & \BM(Q) }$$
\end{itemize}

The empty sequence provides the unit of $\BM(P)$, even when $P$ is 
not unitary.
Note that $\BM$ is essentially surjective: for any monoid $M$, one
has $\BM(M)\simeq M$.

For any pre-monoid $P$, the monoid $\BM(P)$ can be described by the monoid
presentation with
$P$ as set of generators, and a relation $pq=r$ for all $p,q,r\in P$ such
that $pq=r$ in $P$.

Formally, an element $m\in \BM(P)$ is an equivalence class of sequences
of element of $P$, called {\sf $P$-decompositions of $m$}.
A $P$-decomposition
is {\sf reduced} if there it contains no occurrence of the unit of $P$
(if $P$ is not unitary, the condition is empty).

\begin{defi}
Let $P$ be a pre-monoid. Let $m\in \BM(P)$. 
We denote by $E(m)$ the set of reduced $P$-decompositions of $m$.

We denote by $\leq$ the smallest partial order relation on
$E(m)$ such that, for all $(a_1,\dots,a_n) \in E(m)$ and
for all $i$ such that
$a_ia_{i+1}$ is defined in $P$, we have
$$(a_1,\dots,a_{i-1},a_i,a_{i+1},a_{i+2},\dots,a_n)
\leq (a_1,\dots,a_{i-1},a_ia_{i+1},a_{i+2},\dots,a_n).$$
\end{defi}

There is a classical notion of dimension for posets. Let
$(E,\leq)$ be a poset, let $e_0<\dots < e_n$ be a chain in
$E$; the length of the chain is, by definition, the integer $n$.
The {\sf dimension} of $(E,\leq)$ is set to be the supremum
of the set of lengths of all chains in $E$. This dimension
is an element of $\BZ_{\geq 0} \cup \{\infty \}$. 
It coincides with usual notion
of dimension for the simplicial realization of $E$.

\begin{defi}
Let $P$ be a pre-monoid. We say that $P$ is {\sf atomic}
if and only, for all $p \in P$, $E((p))$ is finite dimensional.
\end{defi}

For monoids, this definition coincides with the usual one.

\subsection{The functor $\BG$}
In a similar way, the embedding functor 
${\mathbf{Grp}} \rightarrow \mathbf{preMon}$ has a left adjoint $\BG$,
acting on objects as follows:
for any pre-monoid $P$, the group $\BG(P)$ can be described by 
the group presentation with
$P$ as set of generators, and a relation $pq=r$ for all $p,q,r\in P$ such
that $pq=r$ in $P$.

We have $\BG \BM \simeq \BG$.
If $M$ is monoid satisfying Ore's condition, $\BG(M)$ is isomorphic to 
the group of fractions of $M$.

\subsection{Generated groups}
\label{ggroup}

A general way of constructing a pre-monoid is from a
pair $(G,A)$ where $G$ is a group and $A\subset G$ generates
$G$ as a monoid (we call such a pair a {\sf generated group}).
Let $(G,A)$ be a generated group.
An {\sf $A$-decomposition} of $g\in G$ is a sequence
$(a_1,\dots,a_n)\in A^*$ such that $g=a_1\dots a_n$.
An $A$-decomposition of $g$ of minimal length is said to be 
{\sf reduced}. We denote by $\Red_A(g)$ the set of reduced decompositions
of $g$. We denote by $l_A(g)$ the common length of the elements of $\Red_A(g)$.
The function $l_A$ is sub-additive: for all $g,h\in G$, we have
$$l_A(gh) \leq l_A(g) + l_A(h).$$
We write $g\prec_A h$ if $l_A(g) + l_A(g^{-1}h)=l_A(h)$,
and $g \succ_A h$ if $l_A(gh^{-1}) + l_A(h) = l_A(g)$.

\begin{defi}
Let $(G,A)$ be a generated group.
Let $g\in G$. 
We say that $g$ is {\sf $A$-balanced} (or simply {\sf balanced}) if
$\forall h\in G, h\prec_A g \Leftrightarrow g \succ_A h$.
\end{defi}

Let $g$ be a balanced element of $G$. The set
$$\{h \in G | h\prec_A g\}=\{h \in G | g \succ_A h\}$$
is denoted by $ P_{G,A,g}$ (or simply by $P_g$).
Let $$D_g:=\{(h,h')\in P_g\times P_g | l_A(hh') = l_A(h) + l_A(h') \}$$ and
let $m_g$ be the restriction of the group product to $D_g$.
The triple $(P_g,D_g,m_g)$ is a unitary pre-monoid
(to prove the associativity axiom, use the
fact that $g$ is balanced).

\begin{defi}
The pre-monoid $(P_g,D_g,m_g)$ (or simply $P_g$) is called
{\sf pre-monoid of divisors of $g$} in $(G,A)$.
\end{defi}

Note that the restriction to $P_g$ of the
relation $\prec_A$ (resp. $\succ_A$) is really the left (resp. right)
divisibility relation for the pre-monoid structure.

\begin{defi}
A pre-monoid $P$ is said to be {\sf $\BM$-cancellative} if
$$\forall m \in \BM(P), \forall p,q\in P, 
((pm=qm) \; \text{or} \; (mp=mq)) \Rightarrow p=q.$$
\end{defi}

Note that this is formally weaker than the cancellativity of $\BM(P)$.
A first property of divisors pre-monoids is:

\begin{lemma}
\label{pregar}
Let $(G,A)$ be a presented group.  Let $g$ be a balanced element
of $G$. The pre-monoid $P_g$ is
$\BM$-cancellative.
\end{lemma}

\begin{proof}
Since the defining relations are valid in $G$, the
monoid $\BM(P_g)$ comes equipped with a natural 
morphism $\pi: \BM(P_g) \rightarrow G$.
If for example $pm=qm$, then $\pi(p)\pi(m)=\pi(q)\pi(m)$ in $G$.
Since $G$ is cancellative, $\pi(p)=\pi(q)$. To obtain the first
claim, observe
that the composition of the natural pre-monoid
morphism $P_g \rightarrow \BM(P_g)$ with $\pi$ is the restriction of
identity map of $G$.
\end{proof}

\subsection{Garside monoids} The terminology has been fluctuating in the
recent years,
between several non-equivalent but similar sets of axioms. The
following version seems to emerge as ``consensual''.

\begin{defi}
A monoid $M$ is a {\sf Garside monoid} if:
\begin{itemize}
\item the monoid $M$ is atomic;
\item the monoid $M$ is left and right cancellative;
\item the posets $(M,\prec)$ and $(M,\succ)$ are lattices;
\item there exists an element $\Delta\in M$ such
$$\forall m \in M, (m \prec \Delta ) \Leftrightarrow (\Delta \succ m),$$
and $\{ m \in M | m \prec \Delta\}$ is finite and generates $M$.
(An element $\Delta$ satisfying this property is called 
a {\sf Garside element}.)
\end{itemize}
\end{defi}

Saying that $(M,\prec)$ and $(M,\succ)$ are lattices can be rephrased,
in arithmetical terms, as the existence of left and right
lcm's and gcd's.

Let $M$ be a Garside monoid, with Garside element $\Delta$.
Let $P:=\{m\in M | m\prec \Delta\}$. View $P$ as a pre-monoid,
the product of $p,q\in P$ being defined as the product $pq$ in $M$
(when $pq\in P$; otherwise, the product is not defined).
We call a pre-monoid $P$ obtained this way a {\sf Garside pre-monoid}.
The monoid $M$ can be recovered from $P$: we have $M\simeq \BM(P)$.

In \cite{BDM} is given an axiomatic characterization of Garside
pre-monoids (axioms (i) -- (vi) + existence of a common multiple).
As J. Michel pointed to us, in the context of 
generated groups, most of them are straightforward:

\begin{theo}
\label{BDM}
Let $(G,A)$ be a finite generated group.
Let $g$ be a balanced element in $G$. Assume that $A\subset P_g$,
and that all pairs $a,b\in A$ have a left lcm and a right lcm in $A$.
Then $P_g$ is a Garside pre-monoid.
As a consequence, $\BM(P_g)$ is a Garside monoid.
\end{theo}

The existence of left/right lcm's for pairs of elements
of $A$ follows, for example, if $(P_g,\prec_A)$ and $(P_g,\succ_A)$
are lattices. Conversely, a consequence of the theorem is that
if pairs of elements
of $A$ have left/right lcm's, then $(P_g,\prec_A)$ and $(P_g,\succ_A)$
are lattices. 

This theorem is a convenient tool,
hiding most of the technical machinery (the long list of easy axioms).
But the whole issue remains to check
that $(P_g,\prec_A)$ and $(P_g, \succ_A)$ are lattices.

\begin{proof}
The pre-monoid $P_g$ is unitary (axioms (i) and (ii) of \cite{BDM});
the length function $l_A$ satisfies axiom (iii).
With the assumption $A\subset P_g$, $A$ is the
set of atoms of $P_g$.
The existence of left and right lcm's for elements of $A$ 
are exactly axioms (iv) and (iv').
Axiom (v): let $h\in P_g$, $a,b\in A$,
such that $h\prec_A g$, $ha\prec_A g$ and $hb\prec_A g$;
then $a\prec_A h^{-1}g$, $b\prec_A h^{-1}g$, so $\text{lcm}(a,b)
\prec_A h^{-1}g$ and $h\text{lcm}(a,b)\prec_Ag$.
Axiom (vi) is $\BM$-cancellativity, which we have proved in Lemma
\ref{pregar}.
The element $g$ is a common multiple of all elements of $P_g$.
We conclude using Theorem 2.24 in \cite{BDM}.
\end{proof}

It would be interesting to characterize Garside monoids arising
from triples $(G,A,g)$.

{\bf \flushleft Basic example: the classical braid monoid.}
Let $(W,S)$ be a finite Coxeter system; we view it as a generated
group.
Some crucial results
from \cite{deligne} and \cite{brieskornsaito} show that
the longest element $w_0$ is $S$-balanced (actually, $P_{w_0}= W$
as sets), and that
the posets $(P_{w_0},\prec_S)$ and $(P_{w_0},\succ_S)$ are lattices.
We have $\BB_+(W,S) \simeq \BM(P_{w_0})$.
Our construction of the dual monoid is similar, $S$ being replaced
by $T$ and $w_0$ by a Coxeter element $c$.

\subsection{Properties of Garside monoids}
Let us conclude this section by compiling some remarkable properties
of Garside monoids. 
Any Garside monoid $M$ satisfies the {\sf embedding property}, \ie,
the canonical map $M\rightarrow \BG(M)$ is injective. This implies
that $M$ is cancellative.
Any finite subset of $M$ admits a right lcm, a left lcm, a left gcd
and a right gcd. In particular, $M$ satisfies Ore's conditions on the
left and on the right.
In all examples considered here, the lcm of the atoms is a Garside element.
Let $\Delta$ be a Garside element, with set of divisors $P$.
Any element $m\in M$ has a unique decomposition as a product
$m=p_1\dots p_k$ of elements of $P$ such that, for all 
$i\in\{ 1,\dots,k\}$, $p_i$
is the left gcd of $\Delta$ and $p_i\dots p_k$. The sequence $(p_1,\dots,p_k)$
is called the {\sf normal form} of $m$. One has a similar notion
in $\BG(P)$. This gives rise to solutions of the word problem.
A sequence $(p_1,\dots,p_k)$ is the normal form of $p_1\dots p_k$ if and
only if, for all 
$i\in\{ 1,\dots,k-1\}$, $p_i$
is the left gcd of $\Delta$ and $p_ip_{i+1}$. In other words, the normality
can be checked locally, by looking at consecutive terms. This has important
algorithmic consequences ($\BG(P)$ is biautomatic).
The conjugation action by $\Delta$ on $\BG(M)$ restricts to an automorphism
of $P$. In particular, it is of finite order $d$. We call it the {\sf diagram
automorphism}, by analogy with the case of the classical braid
monoid. It is easy to describe the submonoid of fixed points under 
a given power of
the diagram automorphism.
The element $\Delta^d$ is central in $\BG(M)$.
Some other properties are given in section \ref{applications} of this
paper.

\end{document}